# Penalized nonparametric mean square estimation of the coefficients of diffusion processes

FABIENNE COMTE[*], VALENTINE GENON-CATALOT[**] and YVES ROZENHOLC[†]

*Laboratoire MAP5 (CNRS UMR 8145), UFR de Mathématiques et Informatique, Université Paris 5 René Descartes, 45 rue des Saints Pères, 75270 Paris Cedex 06, France.*
*E-mail:* [*]*fabienne.comte@univ-paris5.fr;* [**]*Valentine.Genon-Catalot@math-info.univ-paris5.fr;*
[†]*yves.rozenholc@math-info.univ-paris5.fr*

We consider a one-dimensional diffusion process $(X_t)$ which is observed at $n+1$ discrete times with regular sampling interval $\Delta$. Assuming that $(X_t)$ is strictly stationary, we propose non-parametric estimators of the drift and diffusion coefficients obtained by a penalized least squares approach. Our estimators belong to a finite-dimensional function space whose dimension is selected by a data-driven method. We provide non-asymptotic risk bounds for the estimators. When the sampling interval tends to zero while the number of observations and the length of the observation time interval tend to infinity, we show that our estimators reach the minimax optimal rates of convergence. Numerical results based on exact simulations of diffusion processes are given for several examples of models and illustrate the qualities of our estimation algorithms.

## 1. Introduction

In this paper, we consider the following problem. Let $(X_t)_{t\geq 0}$ be a one-dimensional diffusion process with dynamics described by the stochastic differential equation:

$$\mathrm{d}X_t = b(X_t)\,\mathrm{d}t + \sigma(X_t)\,\mathrm{d}W_t, \qquad t \geq 0, \quad X_0 = \eta, \tag{1}$$

where $(W_t)$ is a standard Brownian motion and $\eta$ is a random variable independent of $(W_t)$. Assuming that the process is strictly stationary (and ergodic), and that a discrete observation $(X_{k\Delta})_{1\leq k\leq n+1}$ of the sample path is available, we wish to construct non-parametric estimators of the drift function $b$ and the (square of the) diffusion coefficient $\sigma^2$.

Our aim is twofold: to construct estimators that have optimal asymptotic properties and that can be implemented through feasible algorithms. Our asymptotic framework is such that the sampling interval $\Delta = \Delta_n$ tends to zero while $n\Delta_n$ tends to infinity as $n$ tends to infinity. Nevertheless, the risk bounds obtained below are non-asymptotic in the sense that they are explicitly given as functions of $\Delta$ or $1/(n\Delta)$ and fixed constants.

Nonparametric estimation of the coefficients of diffusion processes has been widely investigated in recent decades. The first estimators proposed and studied were based on a continuous time observation of the sample path. Asymptotic results were given for ergodic models as the length of the observation time interval tends to infinity: see, for instance, the reference paper by Banon [2], followed by more recent works by Prakasa Rao [30], Spokoiny [31], Kutoyants [28] or Dalalyan [18].

Then discrete sampling of observations was considered, with different asymptotic frameworks, implying different statistical strategies. It is now classical to distinguish between low-frequency and high-frequency data. In the former case, observations are taken at regularly spaced instants with fixed sampling interval $\Delta$ and the asymptotic framework is that the number of observations tends to infinity. Only ergodic models are usually considered. Parametric estimation in this context was studied by Bibby and Sørensen [11], Kessler and Sørensen [27]; see also Bibby *et al.* [12]. A nonparametric approach using spectral methods was investigated in Gobet *et al.* [24], where non-standard nonparametric rates were exhibited.

In high-frequency data, the sampling interval $\Delta = \Delta_n$ between two successive observations is assumed to tend to zero as the number of observations $n$ tends to infinity. Taking $\Delta_n = 1/n$, so that the length of the observation time interval $n\Delta_n = 1$ is fixed, can only lead to estimating the diffusion coefficient consistently. This was done by Hoffmann [25] who generalized results by Jacod [26], Florens-Zmirou [21] and Genon-Catalot *et al.* [22].

Now, estimating both drift and diffusion coefficients requires that the sampling interval $\Delta_n$ tends to zero while $n\Delta_n$ tends to infinity. For ergodic diffusion models, Hoffmann [25] proposes nonparametric estimators using projections on wavelet bases together with adaptive procedures. He exhibits minimax rates and shows that his estimators automatically reach these optimal rates up to logarithmic factors. Hoffmann's estimators are based on computations of some random times which make them difficult to implement.

In this paper, we propose simple nonparametric estimators based on a penalized mean square approach. The method is investigated in detail in Comte and Rozenholc [16, 17] for regression models. We adapt it here to the case of discretized diffusion models. The estimators are chosen to belong to finite-dimensional spaces that include trigonometric, wavelet-generated and piecewise polynomial spaces. The space dimension is chosen by a data-driven method using a penalization device. Due to the construction of our estimators, we measure the risk of an estimator $\hat{f}$ of $f$ (with $f = b, \sigma^2$) by $\mathbb{E}(\|\hat{f} - f\|_n^2)$, where $\|\hat{f} - f\|_n^2 = n^{-1}\sum_{k=1}^n (\hat{f}(X_{k\Delta}) - f(X_{k\Delta}))^2$. We give bounds for this risk (see Theorems 1 and 2). An examination of these bounds as $\Delta = \Delta_n \to 0$ and $n\Delta_n \to +\infty$ shows that our estimators achieve the optimal nonparametric asymptotic rates obtained in Hoffmann [25] without logarithmic loss (when the unknown functions belong to Besov balls). Then we proceed to numerical implementation on simulated data for several examples of models. We emphasize that our simulation method for diffusion processes is not based on approximations (like Euler schemes). Instead, we use the exact retrospective simulation method described in Beskos *et al.* [10] and Beskos and Roberts [9]. Then we apply the algorithms developed in Comte and Rozenholc [16, 17] for nonparametric estimation using piecewise polynomials. The results are convincing even when some of the theoretical assumptions are not fulfilled.

The paper is organized as follows. In Section 2 we describe our framework (model, assumptions and spaces of approximation). Section 3 is devoted to drift estimation, and Section 4 to diffusion coefficient estimation. In Section 5 we study examples and present numerical simulation results that illustrate the performance of estimators. Section 6 contains proofs. In Section 7 a technical lemma is proved.

## 2. Framework and assumptions

### 2.1. Model assumptions

Let $(X_t)_{t\geq 0}$ be a solution of (1) and assume that $n+1$ observations $X_{k\Delta}$, $k=1,\ldots,n+1$, with sampling interval $\Delta$ are available. Throughout the paper, we assume that $\Delta = \Delta_n$ tends to 0 and $n\Delta_n$ tends to infinity as $n$ tends to infinity. To simplify notation, we write $\Delta$ without the subscript $n$. Nevertheless, when speaking of constants, we mean quantities that depend neither on $n$ nor on $\Delta$. We wish to estimate the drift function $b$ and the diffusion coefficient $\sigma^2$ when $X$ is stationary and geometrically $\beta$-mixing. To this end, we consider the following assumptions:

#### *Assumption 1.*

(i) $b \in C^1(\mathbb{R})$ and there exists $\gamma \geq 0$ such that, for all $x \in \mathbb{R}, |b'(x)| \leq \gamma(1+|x|^\gamma)$.
(ii) There exists $b_0$ such that, for all $x$, $|b(x)| \leq b_0(1+|x|)$.
(iii) There exist $d \geq 0, r > 0$ and $R > 0$ such that, for all $|x| \geq R$, $\mathrm{sgn}(x)b(x) \leq -r|x|^d$.

#### *Assumption 2.*

(i) There exist $\sigma_0^2$ and $\sigma_1^2$ such that, for all $x, 0 < \sigma_0^2 \leq \sigma^2(x) \leq \sigma_1^2$ and there exists $L$ such that, for all $(x,y) \in \mathbb{R}^2, |\sigma(x) - \sigma(y)| \leq L|x-y|^{1/2}$.
(ii) $\sigma \in C^2(\mathbb{R})$ and there exists $\gamma \geq 0$ such that, for all $x \in \mathbb{R}$, $|\sigma'(x)| + |\sigma''(x)| \leq \gamma(1+|x|^\gamma)$.

Under Assumptions 1 and 2, equation (1) has a unique strong solution. Note that Assumption 2(ii) is only used for the estimation of $\sigma^2$ and not for $b$. Elementary computations show that the scale density

$$s(x) = \exp\left\{-2\int_0^x \frac{b(u)}{\sigma^2(u)}\,\mathrm{d}u\right\}$$

satisfies $\int_{-\infty}^{\infty} s(x)\,dx = +\infty = \int^{+\infty} s(x)\,dx$, and the speed density $m(x) = 1/(\sigma^2(x)s(x))$ satisfies $\int_{-\infty}^{+\infty} m(x)\,dx = M < +\infty$. Hence, model (1) admits a unique invariant probability $\pi(x)\,dx$ with $\pi(x) = M^{-1}m(x)$. Now we assume the following:

**Assumption 3.** *$X_0 = \eta$ has distribution $\pi$.*

Under the additional Assumption 3, $(X_t)$ is strictly stationary and ergodic. Moreover, it follows from Proposition 1 in Pardoux and Veretennikov [29] that there exist constants $K > 0$, $\nu > 0$ and $\theta > 0$ such that

$$\mathbb{E}(\exp(\nu|X_0|)) < +\infty \quad \text{and} \quad \beta_X(t) \leq Ke^{-\theta t}, \tag{2}$$

where $\beta_X(t)$ denotes the $\beta$-mixing coefficient of $(X_t)$ and is given by

$$\beta_X(t) = \int_{-\infty}^{+\infty} \pi(x)\,dx \|P_t(x, dx') - \pi(x')\,dx'\|_{\mathrm{TV}}.$$

The norm $\|\cdot\|_{\mathrm{TV}}$ is the total variation norm and $P_t$ denotes the transition probability. In particular, $X_0$ has moments of any (positive) order. Now, Assumption 1(i) ensures that, for all $t \geq 0$, $h > 0$ and $k \geq 1$, there exists $c = c(k, \gamma)$ such that

$$\mathbb{E}\left(\sup_{s \in [t, t+h]} |b(X_s) - b(X_t)|^k | \mathcal{F}_t\right) \leq ch^{k/2}(1 + |X_t|^c),$$

where $\mathcal{F}_t = \sigma(X_s, s \leq t)$; for example, Gloter ([23], Proposition A). Thus, taking expectations, there exists $c'$ such that

$$\mathbb{E}\left(\sup_{s \in [t, t+h]} |b(X_s) - b(X_t)|^k\right) \leq c'h^{k/2}. \tag{3}$$

The functions $b$ and $\sigma^2$ are estimated only on a compact set $A$. For simplicity and without loss of generality, we assume from now on that

$$A = [0, 1]. \tag{4}$$

It follows from Assumptions 1, 2(i) and 3 that the stationary density $\pi$ is bounded from below and above on any compact subset of $\mathbb{R}$, and we denote by $\pi_0$, $\pi_1$ two positive real numbers such that

$$0 < \pi_0 \leq \pi(x) \leq \pi_1 \quad \forall x \in A = [0, 1]. \tag{5}$$

### 2.2. Spaces of approximation: piecewise polynomials

We aim to estimate the functions $b$ and $\sigma^2$ of model (1) on $[0, 1]$ using a data-driven procedure. For that purpose, we consider families of finite-dimensional linear subspaces of $\mathbb{L}^2([0,1])$ and compute for each space an associated least squares estimator. Then an adaptive procedure chooses among the resulting collection of estimators the 'best' one, in a sense that will be specified later, through a penalization device.

Several possible collections of spaces are available and discussed in Section 2.3. Now, to be consistent with the algorithm implemented in Section 5, we focus on a specific collection, namely the collection of dyadic regular piecewise polynomial spaces, henceforth denoted by [DP].

We fix an integer $r \geq 0$. Let $p \geq 0$ also be an integer. On each subinterval $I_j = [(j-1)/2^p, j/2^p]$, $j = 1, \ldots, 2^p$, consider $r+1$ polynomials of degree $\ell$, $\varphi_{j,\ell}(x)$, $\ell = 0, 1, \ldots, r$, and set $\varphi_{j,\ell}(x) = 0$ outside $I_j$. The space $S_m$, $m = (p, r)$, is defined as generated by the $D_m = 2^p(r+1)$ functions $(\varphi_{j,\ell})$. A function $t$ in $S_m$ may be written as

$$t(x) = \sum_{j=1}^{2^p} \sum_{\ell=0}^{r} t_{j,\ell} \varphi_{j,\ell}(x).$$

The collection of spaces $(S_m, m \in \mathcal{M}_n)$ is such that

$$\mathcal{M}_n = \{m = (p, r), p \in \mathbb{N}, r \in \{0, 1, \ldots, r_{\max}\}, 2^p(r_{\max} + 1) \leq N_n\}. \tag{6}$$

In other words, $D_m \leq N_n$, where $N_n \leq n$. The maximal dimension $N_n$ is subject to additional constraints given below. The role of $N_n$ is to bound all dimensions $D_m$, even when $m$ is random. In practice, it corresponds to the maximal number of coefficients to estimate. Thus it must not be too large.

More concretely, consider the orthogonal collection in $\mathbb{L}^2([-1,1])$ of Legendre polynomials $(Q_\ell, \ell \geq 0)$, where the degree of $Q_\ell$ is equal to $\ell$, generating $\mathbb{L}^2([-1,1])$; see Abramowitz and Stegun ([1], page 774). These satisfy $|Q_\ell(x)| \leq 1$, for all $x \in [-1,1]$, $Q_\ell(1) = 1$ and $\int_{-1}^{1} Q_\ell^2(u)\,\mathrm{d}u = 2/(2\ell+1)$. Then we set $P_\ell(x) = (2\ell+1)^{1/2} Q_\ell(2x-1)$ to obtain an orthonormal basis of $\mathbb{L}^2([0,1])$. Finally,

$$\varphi_{j,\ell}(x) = 2^{p/2} P_\ell(2^p x - j + 1)\mathbb{1}_{I_j}(x), \qquad j = 1, \ldots, 2^p, \ell = 0, 1, \ldots, r.$$

The space $S_m$ has dimension $D_m = 2^p(r+1)$, and its orthonormal basis described above satisfies

$$\left\|\sum_{j=1}^{2^p}\sum_{\ell=0}^{r}\varphi_{j,\ell}^2\right\|_\infty \leq D_m(r+1) \leq D_m(r_{\max}+1).$$

Hence, for all $t \in S_m$, $\|t\|_\infty \leq (r_{\max}+1)^{1/2} D_m^{1/2} \|t\|$, where $\|t\|^2 = \int_0^1 t^2(x)\,\mathrm{d}x$, for $t$ in $\mathbb{L}^2([0,1])$, a property which is essential for the proofs.

### 2.3. Other spaces of approximation

From both theoretical and practical points of view, other spaces can be considered, such as the *trigonometric spaces* [T], where $S_m$ is generated by $\{1, 2^{1/2}\cos(2\pi j x), 2^{1/2}\sin(2\pi j x)$ for $j = 1, \ldots, m$}, has dimension $D_m = 2m + 1$ and $m \in \mathcal{M}_n = \{1, \ldots, [n/2] - 1\}$; and the *dyadic wavelet-generated spaces* [W] with regularity $r$ and compact support, as described, for example, in Daubechies [19], Donoho *et al.* [20] or Hoffmann [25].

The key properties that must be fulfilled to fit in our framework are the following:

($\mathcal{H}_1$) *Norm connection*: $(S_m)_{m \in \mathcal{M}_n}$ is a collection of finite-dimensional linear subspaces of $\mathbb{L}^2([0,1])$, with dimension $\dim(S_m) = D_m$ such that $D_m \leq N_n \leq n$, for all $m \in \mathcal{M}_n$, and satisfying:

$$\text{There exists } \Phi_0 > 0 \text{ such that, for all } m \in \mathcal{M}_n, \text{ for all } t \in S_m, \quad \|t\|_\infty \leq \Phi_0 D_m^{1/2} \|t\|. \tag{7}$$

An orthonormal basis of $S_m$ is denoted by $(\varphi_\lambda)_{\lambda \in \Lambda_m}$, where $|\Lambda_m| = D_m$. It follows from Birgé and Massart [13] that property (7) in the context of ($\mathcal{H}_1$) is equivalent to:

$$\text{There exists } \Phi_0 > 0 \text{ such that } \left\| \sum_{\lambda \in \Lambda_m} \varphi_\lambda^2 \right\|_\infty \leq \Phi_0^2 D_m. \tag{8}$$

Thus, for the collection [DP], (8) holds with $\Phi_0^2 = r_{\max} + 1$. Moreover, for results concerning adaptive estimators, we need an additional assumption:

($\mathcal{H}_2$) *Nesting condition*: $(S_m)_{m \in \mathcal{M}_n}$ is a collection of models such that there exists a space denoted by $\mathcal{S}_n$, belonging to the collection, with $S_m \subset \mathcal{S}_n$ for all $m \in \mathcal{M}_n$. We denote by $N_n$ the dimension of $\mathcal{S}_n$: $\dim(\mathcal{S}_n) = N_n$ $(\forall m \in \mathcal{M}_n, D_m \leq N_n)$.

As far as possible below, we keep the notation general to allow extensions to spaces of approximation other than [DP].

## 3. Drift estimation

### 3.1. Drift estimators: non-adaptive case

Let

$$Y_{k\Delta} = \frac{X_{(k+1)\Delta} - X_{k\Delta}}{\Delta} \quad \text{and} \quad Z_{k\Delta} = \frac{1}{\Delta} \int_{k\Delta}^{(k+1)\Delta} \sigma(X_s) \, \mathrm{d}W_s. \tag{9}$$

The following standard regression-type decomposition holds:

$$Y_{k\Delta} = b(X_{k\Delta}) + Z_{k\Delta} + \frac{1}{\Delta} \int_{k\Delta}^{(k+1)\Delta} (b(X_s) - b(X_{k\Delta})) \, \mathrm{d}s,$$

where $b(X_{k\Delta})$ is the main term, $Z_{k\Delta}$ the noise term and the last term is a negligible residual.

Now, for $S_m$ a space of the collection $\mathcal{M}_n$ and for $t \in S_m$, we consider the following regression contrast:

$$\gamma_n(t) = \frac{1}{n} \sum_{k=1}^{n} [Y_{k\Delta} - t(X_{k\Delta})]^2. \tag{10}$$

The estimator belonging to $S_m$ is defined as

$$\hat{b}_m = \arg\min_{t \in S_m} \gamma_n(t). \tag{11}$$

A minimizer of $\gamma_n$ in $S_m$, $\hat{b}_m$ always exists but may not be unique. Indeed, in some common situations the minimization of $\gamma_n$ over $S_m$ leads to an affine space of solutions. Consequently, it becomes impossible to consider a classical $\mathbb{L}^2$-risk for the 'least squares estimator' of $b$ in $S_m$. In contrast, the random $\mathbb{R}^n$-vector $(\hat{b}_m(X_\Delta), \ldots, \hat{b}_m(X_{n\Delta}))'$ is always uniquely defined. Indeed, let us denote by $\Pi_m$ the orthogonal projection (with respect to the inner product of $\mathbb{R}^n$) onto the subspace $\{(t(X_\Delta), \ldots, t(X_{n\Delta}))', t \in S_m\}$ of $\mathbb{R}^n$. Then $(\hat{b}_m(X_\Delta), \ldots, \hat{b}_m(X_{n\Delta}))' = \Pi_m Y$, where $Y = (Y_\Delta, \ldots, Y_{n\Delta})'$. This is the reason why we define the risk of $\hat{b}_m$ by

$$\mathbb{E}\left[\frac{1}{n} \sum_{k=1}^{n} \{\hat{b}_m(X_{k\Delta}) - b(X_{k\Delta})\}^2\right] = \mathbb{E}(\|\hat{b}_m - b\|_n^2),$$

where

$$\|t\|_n^2 = \frac{1}{n} \sum_{k=1}^{n} t^2(X_{k\Delta}). \tag{12}$$

Thus, our risk is the expectation of an empirical norm. Note that, for a deterministic function $t$, $\mathbb{E}(\|t\|_n^2) = \|t\|_\pi^2 = \int t^2(x) \, \mathrm{d}\pi(x)$ where $\pi$ denotes the stationary law. In view of (5), the $\mathbb{L}^2$-norm, $\|\cdot\|$, and the $\mathbb{L}^2(\pi)$-norm, $\|\cdot\|_\pi$, are equivalent for $A$-supported functions.

### 3.2. Risk of the non-adaptive drift estimator

Using (9), (10) and (12), we have

$$\gamma_n(t) - \gamma_n(b) = \|t - b\|_n^2 + \frac{2}{n} \sum_{k=1}^{n} (b-t)(X_{k\Delta}) Z_{k\Delta}$$

$$+ \frac{2}{n\Delta} \sum_{k=1}^{n} (b-t)(X_{k\Delta}) \int_{k\Delta}^{(k+1)\Delta} (b(X_s) - b(X_{k\Delta})) \, \mathrm{d}s.$$

In view of this decomposition, we define the centred empirical process

$$\nu_n(t) = \frac{1}{n} \sum_{k=1}^{n} t(X_{k\Delta}) Z_{k\Delta}. \tag{13}$$

Now denote by $b_m$ the orthogonal projection of $b$ onto $S_m$. By definition of $\hat{b}_m$, $\gamma_n(\hat{b}_m) \leq \gamma_n(b_m)$. So $\gamma_n(\hat{b}_m) - \gamma_n(b) \leq \gamma_n(b_m) - \gamma_n(b)$. This implies

$$\|\hat{b}_m - b\|_n^2 \leq \|b_m - b\|_n^2 + 2\nu_n(\hat{b}_m - b_m)$$
$$+ \frac{2}{n\Delta} \sum_{k=1}^{n} (\hat{b}_m - b_m)(X_{k\Delta}) \int_{k\Delta}^{(k+1)\Delta} (b(X_s) - b(X_{k\Delta}))\,ds.$$

The functions $\hat{b}_m$ and $b_m$ being $A$-supported, we can cancel the terms $\|b\mathbb{1}_{A^c}\|_n^2$ that appear in both sides of the inequality. This yields

$$\|\hat{b}_m - b\mathbb{1}_A\|_n^2 \leq \|b_m - b\mathbb{1}_A\|_n^2 + 2\nu_n(\hat{b}_m - b_m)$$
$$+ \frac{2}{n\Delta} \sum_{k=1}^{n} (\hat{b}_m - b_m)(X_{k\Delta}) \int_{k\Delta}^{(k+1)\Delta} (b(X_s) - b(X_{k\Delta}))\,ds. \tag{14}$$

On the basis of this inequality, we obtain the following result.

**Proposition 1.** *Let $\Delta = \Delta_n$ be such that $\Delta_n \to 0$, $n\Delta_n / \ln^2(n) \to +\infty$ when $n \to +\infty$. Suppose that Assumptions 1, 2(i) and 3 hold and consider a space $S_m$ in the collection [DP] with $N_n = o(n\Delta/\ln^2(n))$ ($N_n$ is defined in $(\mathcal{H}_2)$). Then the estimator $\hat{b}_m$ of $b$ is such that*

$$\mathbb{E}(\|\hat{b}_m - b_A\|_n^2) \leq 7\pi_1 \|b_m - b_A\|^2 + K\frac{\mathbb{E}(\sigma^2(X_0)) D_m}{n\Delta} + K'\Delta + \frac{K''}{n\Delta}, \tag{15}$$

*where $b_A = b\mathbb{1}_A$ and $K, K'$ and $K''$ are positive constants.*

As a consequence, it is natural to select the dimension $D_m$ that leads to the best compromise between the squared bias term $\|b_m - b_A\|^2$ and the variance term of order $D_m/(n\Delta)$.

To compare the result of Proposition 1 with the optimal nonparametric rates exhibited by Hoffmann [25], let us assume that $b_A$ belongs to a ball of some Besov space, $b_A \in \mathcal{B}_{\alpha,2,\infty}([0,1])$, and that $r + 1 \geq \alpha$. Then, for $\|b_A\|_{\alpha,2,\infty} \leq L$, we have $\|b_A - b_m\|^2 \leq C(\alpha, L) D_m^{-2\alpha}$. Thus, choosing $D_m = (n\Delta)^{1/(2\alpha+1)}$, we obtain

$$\mathbb{E}(\|\hat{b}_m - b_A\|_n^2) \leq C(\alpha, L)(n\Delta)^{-2\alpha/(2\alpha+1)} + K'\Delta + \frac{K''}{n\Delta}. \tag{16}$$

The first term $(n\Delta)^{-2\alpha/(2\alpha+1)}$ is exactly the optimal nonparametric rate (see Hoffmann [25]). Moreover, under the standard condition $\Delta = o(1/(n\Delta))$, the last two terms in (15) are $O(1/(n\Delta))$, which is negligible with respect to $(n\Delta)^{-2\alpha/(2\alpha+1)}$.

Proposition 1 holds for the wavelet basis [W] under the same assumptions. For the trigonometric basis [T], the additional constraint $N_n \leq O((n\Delta)^{1/2}/\ln(n))$ is necessary. Hence, when working with these bases, if $b_A \in \mathcal{B}_{\alpha,2,\infty}([0,1])$ as above, the optimal rate is reached for the same choice for $D_m$, under the additional constraint that $\alpha > 1/2$ for [T]. It is worth stressing that $\alpha > 1/2$ automatically holds under Assumption 1.

### 3.3. Adaptive drift estimator

As a second step, we must ensure an automatic selection of $D_m$, which does not use any knowledge of $b$, and in particular which does not require $\alpha$ to be known. The standard selection is

$$\hat{m} = \arg \min_{m \in \mathcal{M}_n} [\gamma_n(\hat{b}_m) + \mathrm{pen}(m)], \qquad (17)$$

with $\mathrm{pen}(m)$ a penalty to be chosen appropriately. We denote by $\hat{b}_{\hat{m}}$ the resulting estimator and we need to determine $\mathrm{pen}(\cdot)$ such that, ideally,

$$\mathbb{E}(\|\hat{b}_{\hat{m}} - b_A\|_n^2) \leq C \inf_{m \in \mathcal{M}_n} \left( \|b_A - b_m\|^2 + \frac{\mathbb{E}(\sigma^2(X_0)) D_m}{n\Delta} \right) + K'\Delta + \frac{K''}{n\Delta},$$

with $C$ a constant which should not be too large. We almost achieve this aim.

**Theorem 1.** *Let $\Delta = \Delta_n$ be such that $\Delta_n \to 0$, $n\Delta_n/\ln^2(n) \to +\infty$ when $n \to +\infty$. Suppose that Assumptions 1, 2(i) and 3 hold and consider the nested collection of models* [DP] *with maximal dimension $N_n = o(n\Delta/\ln^2(n))$. Let*

$$\mathrm{pen}(m) \geq \kappa \sigma_1^2 \frac{D_m}{n\Delta}, \qquad (18)$$

*where $\kappa$ is a universal constant. Then the estimator $\hat{b}_{\hat{m}}$ of $b$ with $\hat{m}$ defined in* (17) *is such that*

$$\mathbb{E}(\|\hat{b}_{\hat{m}} - b_A\|_n^2) \leq C \inf_{m \in \mathcal{M}_n} (\|b_m - b_A\|^2 + \mathrm{pen}(m)) + K'\Delta + \frac{K''}{n\Delta}. \qquad (19)$$

Some comments are in order. It is possible to choose $\mathrm{pen}(m) = \kappa \sigma_1^2 D_m/(n\Delta)$, but this is not what is done in practice. It is better to calibrate additional terms. This is explained in Section 5.2. The constant $\kappa$ in the penalty is numerical and must be calibrated for the problem. Its value is usually adapted by intensive simulation experiments. This point is also discussed in Section 5.2. From (15), one would expect to obtain $\mathbb{E}(\sigma^2(X_0))$ instead of $\sigma_1^2$ in (18): we do not know if this is the consequence of technical problems or if it is a structural result. Another important point is that $\sigma_1^2$ is unknown. In practice, we just replace it by a rough estimator (see Section 5.2).

From (19), we deduce that the adaptive estimator automatically realizes the bias–variance compromise: whenever $b_A$ belongs to some Besov ball (see (16)), if $r + 1 \geq \alpha$ and $n\Delta^2 = o(1)$, $\hat{b}_{\hat{m}}$ achieves the optimal corresponding nonparametric rate, without logarithmic loss, contrary to Hoffmann's adaptive estimator (see Hoffmann [25], page 159, Theorem 5). As mentioned above, Theorem 1 holds for the basis [W] and, if $N_n = o((n\Delta)^{1/2}/\ln(n))$, for [T].

## 4. Adaptive estimation of the diffusion coefficient

### 4.1. Diffusion coefficient estimator: non-adaptive case

To estimate $\sigma^2$ on $A = [0,1]$, we define

$$\hat{\sigma}_m^2 = \arg\min_{t \in S_m} \breve{\gamma}_n(t) \qquad \text{with } \breve{\gamma}_n(t) = \frac{1}{n}\sum_{k=1}^n [U_{k\Delta} - t(X_{k\Delta})]^2, \tag{20}$$

and

$$U_{k\Delta} = \frac{(X_{(k+1)\Delta} - X_{k\Delta})^2}{\Delta}. \tag{21}$$

For diffusion coefficient estimation under our asymptotic framework, it is now well known that rates of convergence are faster than for drift estimation. This is the reason why the regression-type equation has to be more precise than for $b$. Let us set

$$\psi = 2\sigma'\sigma b + [(\sigma')^2 + \sigma\sigma'']\sigma^2. \tag{22}$$

Some computations using Itô's formula and Fubini's theorem lead to

$$U_{k\Delta} = \sigma^2(X_{k\Delta}) + V_{k\Delta} + R_{k\Delta}$$

where $V_{k\Delta} = V_{k\Delta}^{(1)} + V_{k\Delta}^{(2)} + V_{k\Delta}^{(3)}$, with

$$V_{k\Delta}^{(1)} = \frac{1}{\Delta}\left[\left\{\int_{k\Delta}^{(k+1)\Delta} \sigma(X_s)\,dW_s\right\}^2 - \int_{k\Delta}^{(k+1)\Delta} \sigma^2(X_s)\,ds\right]$$

$$V_{k\Delta}^{(2)} = \frac{2}{\Delta}\int_{k\Delta}^{(k+1)\Delta}((k+1)\Delta - s)\sigma'(X_s)\sigma^2(X_s)\,dW_s,$$

$$V_{k\Delta}^{(3)} = 2b(X_{k\Delta})\int_{k\Delta}^{(k+1)\Delta} \sigma(X_s)\,dW_s,$$

and

$$R_{k\Delta} = \frac{1}{\Delta}\left(\int_{k\Delta}^{(k+1)\Delta} b(X_s)\,ds\right)^2 + \frac{2}{\Delta}\int_{k\Delta}^{(k+1)\Delta}(b(X_s) - b(X_{k\Delta}))\,ds \int_{k\Delta}^{(k+1)\Delta} \sigma(X_s)\,dW_s$$

$$+ \frac{1}{\Delta}\int_{k\Delta}^{(k+1)\Delta} [(k+1)\Delta - s]\psi(X_s)\,ds.$$

Obviously, the main noise term in the above decomposition must be $V_{k\Delta}^{(1)}$, as will be proved below.

### 4.2. Risk of the non-adaptive estimator

As for the drift, we write

$$\breve{\gamma}_n(t) - \breve{\gamma}_n(\sigma^2) = \|\sigma^2 - t\|_n^2 + \frac{2}{n}\sum_{k=1}^n (\sigma^2 - t)(X_{k\Delta})V_{k\Delta} + \frac{2}{n}\sum_{k=1}^n (\sigma^2 - t)(X_{k\Delta})R_{k\Delta}.$$

We denote by $\sigma_m^2$ the orthogonal projection of $\sigma^2$ on $S_m$ and define

$$\breve{\nu}_n(t) = \frac{1}{n}\sum_{k=1}^n t(X_{k\Delta})V_{k\Delta}.$$

Again we use the fact that $\breve{\gamma}_n(\hat{\sigma}_m^2) - \breve{\gamma}_n(\sigma^2) \leq \breve{\gamma}_n(\sigma_m^2) - \breve{\gamma}_n(\sigma^2)$ to obtain

$$\|\hat{\sigma}_m^2 - \sigma^2\|_n^2 \leq \|\sigma_m^2 - \sigma^2\|_n^2 + 2\breve{\nu}_n(\hat{\sigma}_m^2 - \sigma_m^2) + \frac{2}{n}\sum_{k=1}^n (\hat{\sigma}_m^2 - \sigma_m^2)(X_{k\Delta})R_{k\Delta}.$$

Analogously to what was done for the drift, we can cancel on both sides the common term $\|\sigma^2 \mathbb{1}_{A^c}\|_n^2$. This yields

$$\|\hat{\sigma}_m^2 - \sigma_A^2\|_n^2 \leq \|\sigma_m^2 - \sigma_A^2\|_n^2 + 2\breve{\nu}_n(\hat{\sigma}_m^2 - \sigma_m^2) + \frac{2}{n}\sum_{k=1}^n (\hat{\sigma}_m^2 - \sigma_m^2)(X_{k\Delta})R_{k\Delta}. \quad (23)$$

We obtain the following result.

**Proposition 2.** *Let $\Delta = \Delta_n$ be such that $\Delta_n \to 0$, $n\Delta_n/\ln^2(n) \to +\infty$ when $n \to +\infty$. Suppose that Assumptions 1–3 hold and consider a model $S_m$ in the collection [DP] with $N_n = o(n\Delta/\ln^2(n))$, where $N_n$ is defined in $(\mathcal{H}_2)$. Then the estimator $\hat{\sigma}_m^2$ of $\sigma^2$ defined by (20) is such that*

$$\mathbb{E}(\|\hat{\sigma}_m^2 - \sigma_A^2\|_n^2) \leq 7\pi_1\|\sigma_m^2 - \sigma_A^2\|^2 + K\frac{\mathbb{E}(\sigma^4(X_0))D_m}{n} + K'\Delta^2 + \frac{K''}{n}, \quad (24)$$

*where $\sigma_A^2 = \sigma^2 \mathbb{1}_A$, and $K$, $K'$, $K''$ are positive constants.*

Let us make some comments on the rates of convergence. If $\sigma_A^2$ belongs to a ball of some Besov space, say $\sigma_A^2 \in \mathcal{B}_{\alpha,2,\infty}([0,1])$, and $\|\sigma_A^2\|_{\alpha,2,\infty} \leq L$, with $r+1 \geq \alpha$, then $\|\sigma_A^2 - \sigma_m^2\|^2 \leq C(\alpha, L)D_m^{-2\alpha}$. Therefore, if we choose $D_m = n^{1/(2\alpha+1)}$, we obtain

$$\mathbb{E}(\|\hat{\sigma}_m^2 - \sigma_A^2\|_n^2) \leq C(\alpha, L)n^{-2\alpha/(2\alpha+1)} + K'\Delta^2 + \frac{K''}{n}. \quad (25)$$

The first term $n^{-2\alpha/(2\alpha+1)}$ is the optimal nonparametric rate proved by Hoffmann [25]. Moreover, under the standard condition $\Delta^2 = o(1/n)$, the last two terms are $O(1/n)$, that is, negligible with respect to $n^{-2\alpha/(2\alpha+1)}$.

### 4.3. Adaptive diffusion coefficient estimator

As previously, the second step is to ensure an automatic selection of $D_m$, which does not use any knowledge on $\sigma^2$. This selection is done by

$$\hat{m} = \arg\min_{m \in \mathcal{M}_n} [\breve{\gamma}_n(\hat{\sigma}_m^2) + \widetilde{\text{pen}}(m)]. \tag{26}$$

We denote by $\hat{\sigma}_{\hat{m}}^2$ the resulting estimator and we need to determine the penalty $\widetilde{\text{pen}}$ as for $b$. For simplicity, we use the same notation $\hat{m}$ in (26) as in (17) although they are different. We can prove the following theorem.

**Theorem 2.** *Let $\Delta = \Delta_n$ be such that $\Delta_n \to 0$, $n\Delta_n / \ln^2(n) \to +\infty$ when $n \to +\infty$. Suppose that Assumptions 1–3 hold. Consider the nested collection of models* [DP] *with maximal dimension $N_n \leq n\Delta / \ln^2(n)$. Let*

$$\widetilde{\text{pen}}(m) \geq \tilde{\kappa}\sigma_1^4 \frac{D_m}{n}, \tag{27}$$

*where $\tilde{\kappa}$ is a universal constant. Then, the estimator $\hat{\sigma}_{\hat{m}}^2$ of $\sigma^2$ with $\hat{m}$ defined by (26) is such that*

$$\mathbb{E}(\|\hat{\sigma}_{\hat{m}}^2 - \sigma_A^2\|_n^2) \leq C \inf_{m \in \mathcal{M}_n} (\|\sigma_m^2 - \sigma_A^2\|^2 + \widetilde{\text{pen}}(m)) + K'\Delta^2 + \frac{K''}{n}. \tag{28}$$

As for the drift, it is possible to choose $\widetilde{\text{pen}}(m) = \tilde{\kappa}\sigma_1^4 D_m/n$, but this is not what is done in practice. Moreover, making such a choice, it follows from (28) that the adaptive estimator automatically realizes the bias–variance compromise. Whenever $\sigma_A^2$ belongs to some Besov ball (see (25)), if $n\Delta^2 = o(1)$ and $r + 1 \geq \alpha$, $\hat{\sigma}_{\hat{m}}^2$ achieves the optimal corresponding nonparametric rate $n^{-2\alpha/(2\alpha+1)}$, without logarithmic loss, contrary to Hoffmann's adaptive estimator (see Hoffmann [25], page 160, Theorem 6). As mentioned for $b$, Proposition 2 and Theorem 2 hold for the basis [W] under the same assumptions on $N_n$. For [T], $N_n = o((n\Delta)^{1/2}/\ln(n))$ is needed.

## 5. Examples and numerical simulation results

In this section, we consider examples of diffusions and implement the estimation algorithms on simulated data. To simulate sample paths of diffusion, we use the retrospective exact simulation algorithms proposed by Beskos *et al.* [10] and Beskos and Roberts [9]. Contrary to the Euler scheme, these algorithms produce exact simulation of diffusions under some assumptions on the drift and diffusion coefficient. Therefore, we choose our examples in order to meet these conditions in addition with our set of assumptions. For the sake of simplicity, we focus on models that can be simulated by the simplest algorithm of Beskos *et al.* [10], which is called EA1. More precisely, consider a diffusion model given by the stochastic differential equation

$$dX_t = b(X_t)\, dt + \sigma(X_t)\, dW_t. \tag{29}$$

We assume that there is a $C^2$ one-to-one mapping $F$ on $\mathbb{R}$ such that $\xi_t = F(X_t)$ satisfies

$$d\xi_t = \alpha(\xi_t)\, dt + dW_t. \tag{30}$$

To produce an exact realization of the random variable $\xi_\Delta$, given that $\xi_0 = x$, the exact algorithm EA1 requires that $\alpha$ be $C^1$, and $\alpha^2 + \alpha'$ be bounded from below and above. Moreover, setting $A(\xi) = \int^\xi \alpha(u)\, du$, the function

$$h(\xi) = \exp(A(\xi) - (\xi - x)^2/2\Delta) \tag{31}$$

must be integrable on $\mathbb{R}$, and an exact realization of a random variable with density proportional to $h$ must be possible. Provided that the process $(\xi_t)$ admits a stationary distribution that it may also be possible to simulate, using the Markov property, the algorithm can therefore produce an exact realization of a discrete sample $(\xi_{k\Delta}, k = 0, 1, \ldots, n+1)$ in the stationary regime. We deduce an exact realization of $(X_{k\Delta} = F^{-1}(\xi_{k\Delta}), k = 0, \ldots, n+1)$.

In all examples, we estimate the drift function $\alpha(\xi)$ and the constant 1 for models like (30) or both the drift $b(x)$ and the diffusion coefficient $\sigma^2(x)$ for models like (29). Let us note that Assumptions 1–3 are fulfilled for all the models $(\xi_t)$ below. For the models $(X_t)$, the ergodicity and the exponential $\beta$-mixing property hold.

## 5.1. Examples of diffusions

### 5.1.1. Family 1

First, we consider (29) with

$$b(x) = -\theta x, \qquad \sigma(x) = c(1 + x^2)^{1/2}. \tag{32}$$

Standard computations of the scale and speed densities show that the model is positive recurrent for $\theta + c^2/2 > 0$. In this case, its stationary distribution has density

$$\pi(x) \propto \frac{1}{(1 + x^2)^{1+\theta/c^2}}.$$

If $X_0 = \eta$ has distribution $\pi(x)\, dx$, then, setting $\nu = 1 + 2\theta/c^2$, $\nu^{1/2}\eta$ has Student distribution $t(\nu)$ which can be easily simulated.

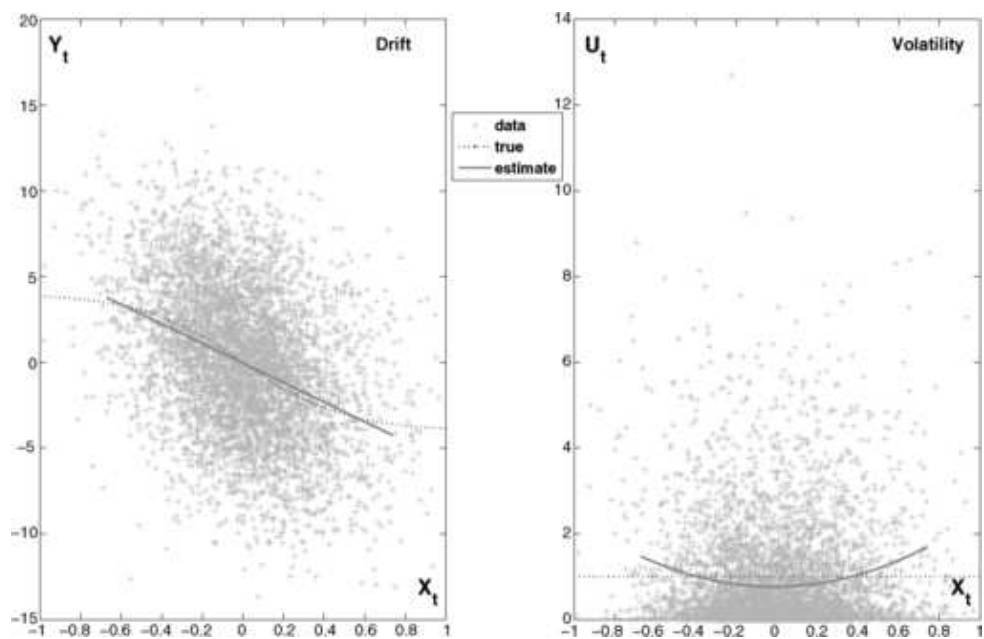

**Figure 1.** $dX_t = -(\theta/c + c/2)\tanh(cX_t) + dW_t$, $n = 5000$, $\Delta = 1/20$, $\theta = 6$, $c = 2$. Dotted line, true; solid line, estimate. The algorithm selects $(p, r)$ equal to $(0, 1)$ for the drift, $(0, 2)$ for $\sigma^2$.

We now consider $F_1(x) = \int_0^x 1/(c(1+x^2)^{1/2})\,dx = \arg\sinh(x)/c$. By the Itô formula, $\xi_t = F_1(X_t)$ satisfies (30) with

$$\alpha(\xi) = -(\theta/c + c/2)\tanh(c\xi). \tag{33}$$

Assumptions 1–3 hold for $(\xi_t)$ with $\xi_0 = F_1(X_0)$. Moreover,

$$\alpha^2(\xi) + \alpha'(\xi) = \{(\theta/c + c/2)^2 + \theta + c^2/2\}\tanh^2(c\xi) - (\theta + c^2/2)$$

is bounded from below and above. And

$$A(\xi) = \int_0^\xi \alpha(u)\,du = -(1/2 + \theta/c^2)\log(\cosh(c\xi)) \leq 0,$$

so that $\exp(A(\xi)) \leq 1$. Therefore, function (31) is integrable for all $x$ and, by a simple rejection method, we can produce a realization of a random variable with density proportional to $h(\xi)$ using a random variable with density $\mathcal{N}(x, \Delta)$.

Note that model (29) satisfies Assumptions 1–3 except that $\sigma^2(x)$ is not bounded from above. Nevertheless, since $X_t = F_1^{-1}(\xi_t) = \sinh(c\xi_t)$, the process $(X_t)$ is exponentially $\beta$-mixing. The upper bound $\sigma_1^2$ that appears explicitly in the penalty function must be replaced by an estimated upper bound.

*5.1.2. Family 2*

For the second family of models, we start with an equation of type (30) where the drift is now (see Barndorff-Nielsen [7])

$$\alpha(\xi) = -\theta \frac{\xi}{(1 + c^2\xi^2)^{1/2}}. \tag{34}$$

The model for $(\xi_t)$ is positive recurrent on $\mathbb{R}$ for $\theta > 0$. Its stationary distribution is given by

$$\pi(\xi)\,\mathrm{d}\xi \propto \exp\left(-2\frac{\theta}{c^2}(1 + c^2\xi^2)^{1/2}\right) = \exp\left(-2\frac{\theta|\xi|}{c}\right)\exp(\varphi(\xi)),$$

where $\exp\varphi(\xi) \leq 1$ so that a random variable with distribution $\pi(\xi)\,\mathrm{d}\xi$ can be simulated by simple rejection method using a double exponential variable with distribution proportional to $\exp(-2\theta|\xi|/c)$. The conditions required to perform an exact simulation of $(\xi_t)$ hold. More precisely, $\alpha^2 + \alpha'$ is bounded from below and above and $A(\xi) = \int_0^\xi \alpha(u)\,\mathrm{d}u = -(\theta/c^2)(1 + c^2\xi^2)^{1/2}$. Hence $\exp(A(\xi)) \leq 1$, (31) is integrable and we can produce a realization of a random variable with density proportional to (31). Lastly, Assumptions 1–3 also hold for this model.

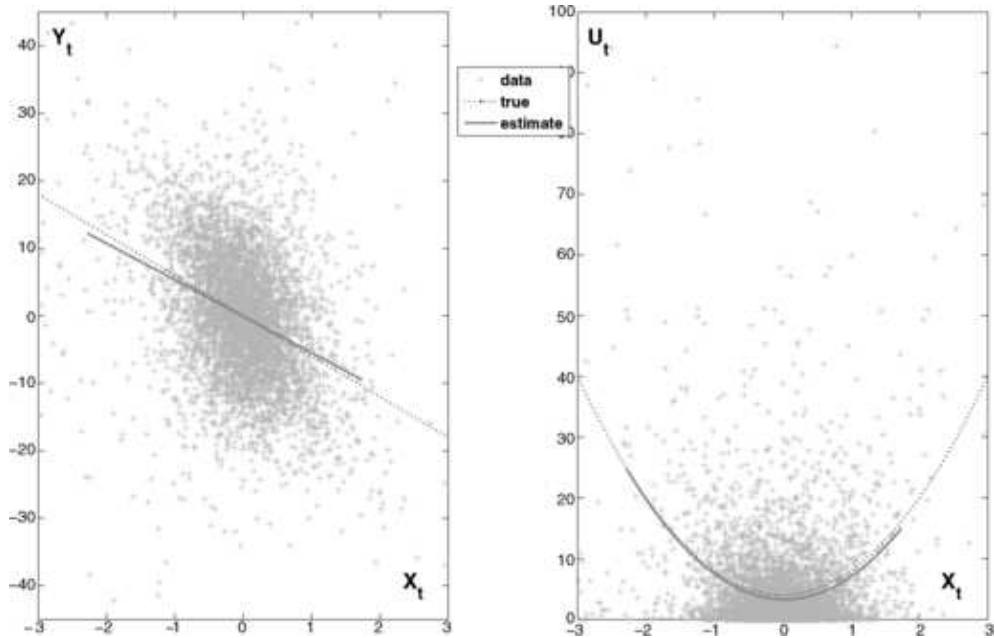

**Figure 2.** $\mathrm{d}X_t = -\theta X_t\,\mathrm{d}t + c(1 + X_t^2)^{1/2}\,\mathrm{d}W_t$, $n = 5000$, $\Delta = 1/20$, $\theta = 6$, $c = 2$. Dotted line, true; solid line, estimate. The algorithm selects $(p, r)$ equal to $(0, 1)$ for the drift, $(0, 2)$ for $\sigma^2$.

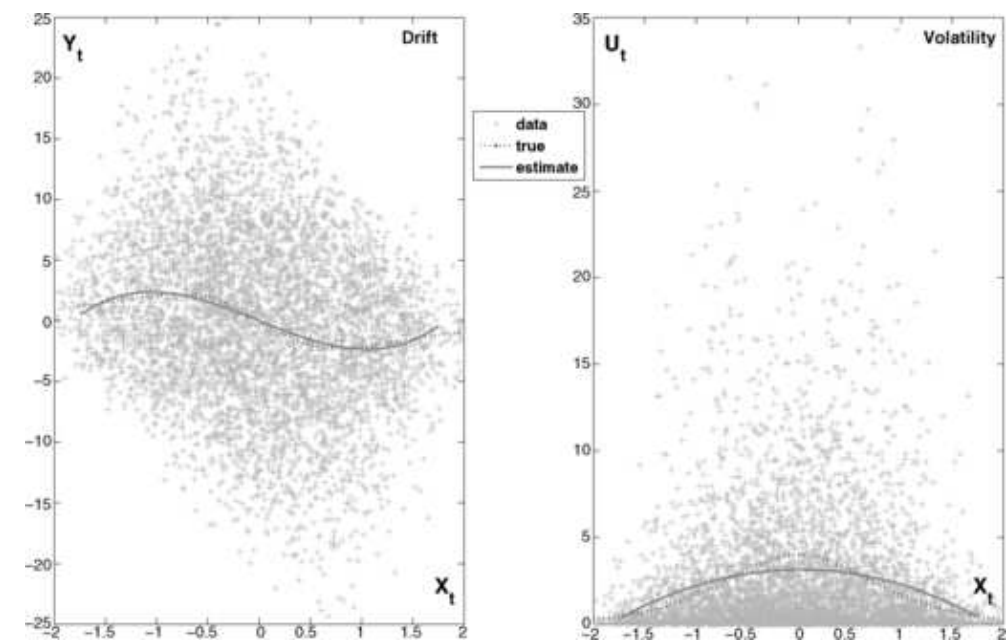

**Figure 3.** $dX_t = -[\theta + c^2/(2\cosh(X_t))](\sinh(X_t)/\cosh^2(X_t))\,dt + (c/\cosh(X_t))\,dW_t$, $n = 5000$, $\Delta = 1/20$, $\theta = 3$, $c = 2$. Dotted line, true; solid line, estimate. The algorithm selects $(p, r)$ equal to $(0, 2)$ for the drift, $(0, 3)$ for $\sigma^2$.

We now consider $X_t = F_2(\xi_t) = \arg\sinh(c\xi_t)$, which satisfies a stochastic differential equation with coefficients

$$b(x) = -\left(\theta + \frac{c^2}{2\cosh(x)}\right)\frac{\sinh(x)}{\cosh^2(x)}, \qquad \sigma(x) = \frac{c}{\cosh(x)}. \tag{35}$$

The process $(X_t)$ is exponentially $\beta$-mixing as $(\xi_t)$. The diffusion coefficient $\sigma(x)$ is not bounded from below but has an upper bound.

To obtain a different shape for the diffusion coefficient, showing two bumps, we consider $X_t = G(\xi_t) = \arg\sinh(\xi_t - 5) + \arg\sinh(\xi_t + 5)$ where $(\xi_t)$ is as in (30)–(34). The function $G(\cdot)$ is invertible and its inverse has the explicit expression

$$G^{-1}(x) = \frac{1}{2^{1/2}\sinh(x)}[49\sinh^2(x) + 100 + \cosh(x)(\sinh^2(x) - 100)]^{1/2}.$$

The diffusion coefficient of $(X_t)$ is given by

$$\sigma(x) = \frac{1}{(1 + (G^{-1}(x) - 5)^2)^{1/2}} + \frac{1}{(1 + (G^{-1}(x) + 5)^2)^{1/2}}. \tag{36}$$

The drift is given by $b(x) = G'(G^{-1}(x))\alpha(G^{-1}(x)) + \frac{1}{2}G''(G^{-1}(x))$.

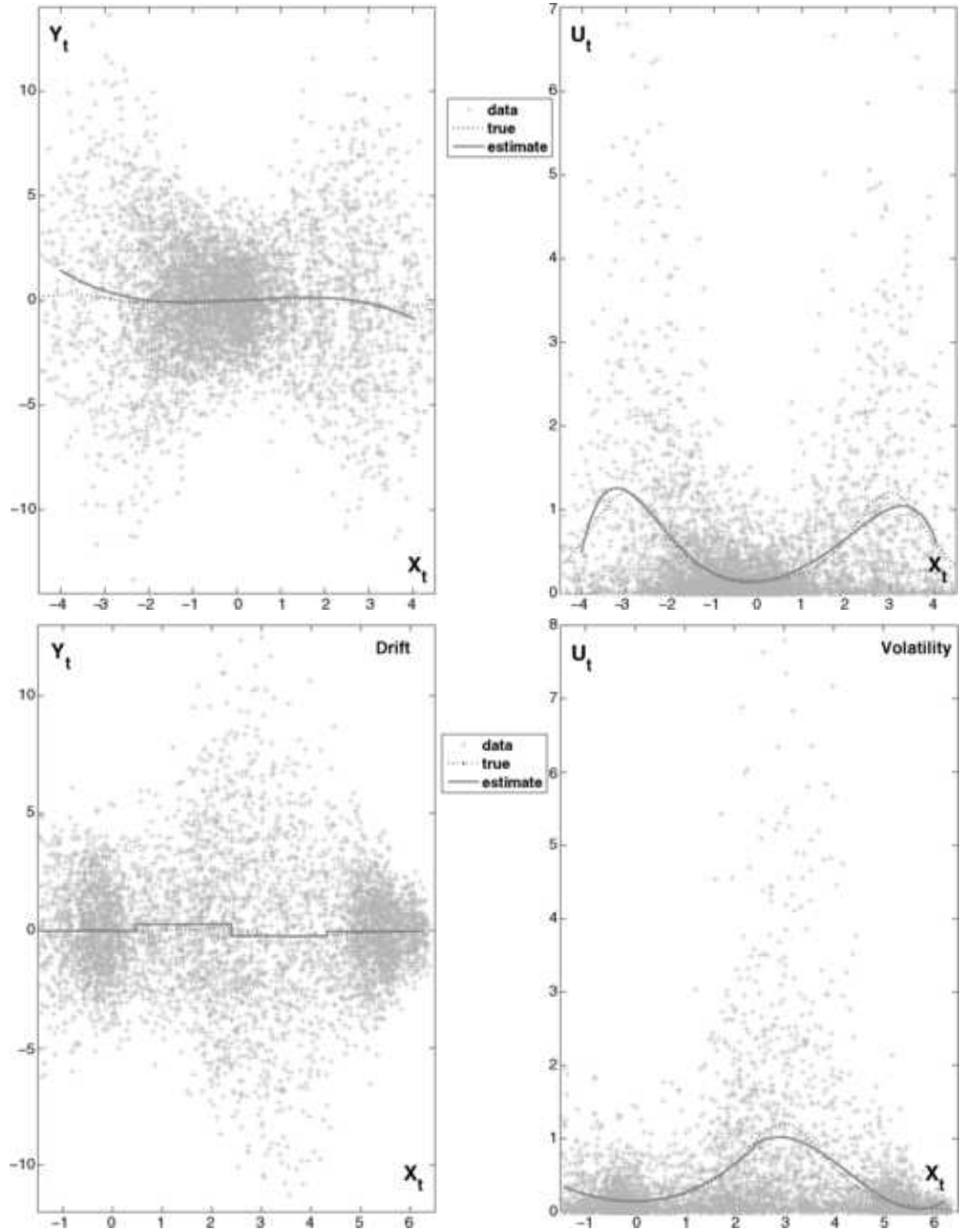

**Figure 4.** Two paths for the two-bumps diffusion coefficient model $X_t = G(\xi_t)$, $\mathrm{d}\xi_t = -\theta\xi_t/(1+c^2\xi_t^2)^{1/2}\,\mathrm{d}t + \mathrm{d}W_t$, $G(x) = \arg\sinh(x-5) + \arg\sinh(x+5)$, $n = 5000$, $\Delta = 1/20$, $\theta = 1$, $c = 10$. Dotted line, true; solid line, estimate. The algorithm selects $(p,r)$ equal to $(0,3)$ (above) and $(2,0)$ (below) for the drift, $(0,6)$ (above) and $(1,3)$ (below) for $\sigma^2$.

## 5.2. Estimation algorithms and numerical results

We do not give here a complete Monte Carlo study but we illustrate how the algorithm works and what kind of estimate it delivers visually.

We consider the regular collection [DP] (see Section 2.2). The algorithm minimizes the mean square contrast and selects the space of approximation in the sense that it selects $p$ and $r$ for integers $p$ and $r$ such that $2^p(r+1) \leq N_n \leq n\Delta/\ln^2(n)$ and $r \in \{0, 1, \ldots, r_{\max}\}$. Note that the degree is global in the sense that it is the same on all the intervals of the subdivision. We take $r_{\max} = 9$ in practice. Moreover, additive (but negligible) correcting terms are classically involved in the penalty (see Comte and Rozenholc [17]). Such terms avoid underpenalization and are in accordance with the fact that the theorems provide lower bounds for the penalty. The correcting terms are asymptotically negligible so they do not affect the rate of convergence. Thus, both penalties contain additional logarithmic terms which have been calibrated in other contexts by intensive simulation experiments (see Comte and Rozenholc [16, 17]).

The constant $\kappa$ in both penalties $\mathrm{pen}(m)$ and $\widetilde{\mathrm{pen}}(m)$ has been set equal to 4.

We retain the idea that the adequate term in the penalty was $\mathbb{E}(\sigma^2(X_0))/\Delta$ for $b$ and $\mathbb{E}(\sigma^4(X_0))$ for $\sigma^2$, instead of those obtained ($\sigma_1^2/\Delta$ and $\sigma_1^4$, respectively). Indeed, in classical regression models, the corresponding coefficient is the variance of the noise. This variance is usually unknown and replaced by a rough estimate. Therefore, in penalties, $\sigma_1^2/\Delta$ and $\sigma_1^4$ are replaced by empirical variances computed using initial estimators $\hat{b}$, $\hat{\sigma}^2$ chosen in the collection and corresponding to a space with medium dimension: $\sigma_1^2/\Delta$ for $\mathrm{pen}(\cdot)$ is replaced $\hat{s}_1^2 = \gamma_n(\hat{b})$ (see (10)); and $\sigma_1^4$ for the other penalty is replaced by $\hat{s}_2^2 = \breve{\gamma}_n(\hat{\sigma}^2)$ (see (20)).

Finally, for $m = (p, r)$, the penalties $\mathrm{pen}(m)$ for $i = 1$ and $\widetilde{\mathrm{pen}}(m)$ for $i = 2$ are given by

$$4\frac{\hat{s}_i^2}{n}2^p(r+1+\ln^{2.5}(r+1)).$$

Figures 1–4 illustrate our simulation results. We have plotted the data points $(X_{k\Delta}, Y_{k\Delta})$ (see (9)) and $(X_{k\Delta}, U_{k\Delta})$ (see (21)), the true functions $b$ and $\sigma^2$ and the estimated functions based on 95% of data points. Parameters have been chosen in the admissible range of ergodicity. The sample size $n = 5000$ and the step size $\Delta = 1/20$ are in accordance with the asymptotic context (large $n$ and small $\Delta$) and may be relevant for applications in finance. It is clear that the estimated functions correspond very well to the true ones.

The simulation of sample paths does not rely on Euler schemes as in the estimation method. Therefore, the data simulation method is disconnected with the estimation procedures and cannot be suspected of being favourable to our estimation algorithm.

## 6. Proofs

### 6.1. Proof of Proposition 1

We recall that for $A$-supported functions, $\|t\|_\pi^2 = \int_A t^2(x)\pi(x)\,\mathrm{d}x$. Starting from (13)–(14), we obtain

$$\|\hat{b}_m - b_A\|_n^2 \leq \|b_m - b_A\|_n^2 + 2\|\hat{b}_m - b_m\|_\pi \sup_{t\in S_m, \|t\|_\pi=1} |\nu_n(t)|$$

$$+ 2\|\hat{b}_m - b_m\|_n \left[\frac{1}{n\Delta^2}\sum_{k=1}^n \left\{\int_{k\Delta}^{(k+1)\Delta}(b(X_s) - b(X_{k\Delta}))\,\mathrm{d}s\right\}^2\right]^{1/2}$$

$$\leq \|b_m - b_A\|_n^2 + \frac{1}{8}\|\hat{b}_m - b_m\|_\pi^2 + 8\sup_{t\in S_m, \|t\|_\pi=1}[\nu_n(t)]^2$$

$$+ \frac{1}{8}\|\hat{b}_m - b_m\|_n^2 + \frac{8}{n\Delta^2}\sum_{k=1}^n\left(\int_{k\Delta}^{(k+1)\Delta}(b(X_s) - b(X_{k\Delta}))\,\mathrm{d}s\right)^2.$$

Because the $\mathbb{L}^2$-norm, $\|\cdot\|_\pi$, and the empirical norm (12) are not equivalent, we must introduce a set on which they are and then prove that this set has small probability. Let us define (see (6))

$$\Omega_n = \left\{\omega \Big/ \left|\frac{\|t\|_n^2}{\|t\|_\pi^2} - 1\right| \leq \frac{1}{2}, \forall t \in \bigcup_{m,m'\in\mathcal{M}_n}(S_m + S_{m'}) \setminus \{0\}\right\}. \tag{37}$$

On $\Omega_n$, $\|\hat{b}_m - b_m\|_\pi^2 \leq 2\|\hat{b}_m - b_m\|_n^2$ and $\|\hat{b}_m - b_m\|_n^2 \leq 2(\|\hat{b}_m - b_A\|_n^2 + \|b_m - b_A\|_n^2)$. Hence, some elementary computations yield:

$$\frac{1}{4}\|\hat{b}_m - b_A\|_n^2 \mathbb{1}_{\Omega_n}$$

$$\leq \frac{7}{4}\|b_m - b_A\|_n^2 + 8\sup_{t\in S_m, \|t\|_\pi=1}[\nu_n(t)]^2 + \frac{8}{n\Delta^2}\sum_{k=1}^n\left(\int_{k\Delta}^{(k+1)\Delta}(b(X_s) - b(X_{k\Delta}))\,\mathrm{d}s\right)^2.$$

Now, using (3), we obtain

$$\mathbb{E}\left(\int_{k\Delta}^{(k+1)\Delta}(b(X_s) - b(X_{k\Delta}))\,\mathrm{d}s\right)^2 \leq \Delta\int_{k\Delta}^{(k+1)\Delta}\mathbb{E}[(b(X_s) - b(X_{k\Delta}))^2]\,\mathrm{d}s \leq c'\Delta^3.$$

Consequently,

$$\mathbb{E}(\|\hat{b}_m - b_A\|_n^2 \mathbb{1}_{\Omega_n}) \leq 7\|b_m - b_A\|_\pi^2 + 32\mathbb{E}\left(\sup_{t\in S_m, \|t\|_\pi=1}[\nu_n(t)]^2\right) + 32c'\Delta. \tag{38}$$

Next, using (5), (7)–(9) and (13), it is easy to see that, since $\|t\|_\pi = 1 \Rightarrow \|t\|^2 \leq 1/\pi_0$,

$$\mathbb{E}\left(\sup_{t \in S_m, \|t\|_\pi = 1}[\nu_n(t)]^2\right) \leq \frac{1}{\pi_0}\mathbb{E}\left(\sup_{t \in S_m, \|t\| \leq 1}[\nu_n(t)]^2\right) \leq \frac{1}{\pi_0}\sum_{\lambda \in \Lambda_m}\mathbb{E}[\nu_n^2(\varphi_\lambda)]$$

$$= \frac{1}{\pi_0 n^2 \Delta^2}\sum_{k=1}^n \mathbb{E}\left\{\sum_{\lambda \in \Lambda_m}\varphi_\lambda^2(X_{k\Delta})\int_{k\Delta}^{(k+1)\Delta}\sigma^2(X_s)\,\mathrm{d}s\right\}$$

$$\leq \frac{\Phi_0^2 D_m}{\pi_0 n^2 \Delta^2}\sum_{k=1}^n \mathbb{E}\left\{\int_{k\Delta}^{(k+1)\Delta}\sigma^2(X_s)\,\mathrm{d}s\right\}$$

$$= \frac{\Phi_0^2 D_m}{\pi_0 n \Delta^2}\mathbb{E}\left(\int_0^\Delta \sigma^2(X_s)\,\mathrm{d}s\right) = \frac{\Phi_0^2 \mathbb{E}(\sigma^2(X_0))D_m}{\pi_0 n \Delta}.$$

Gathering bounds, and using the upper bound $\pi_1$ defined in (5), we obtain

$$\mathbb{E}(\|\hat{b}_m - b_A\|_n^2 \mathbb{1}_{\Omega_n}) \leq 7\pi_1 \|b_m - b_A\|^2 + 32\frac{\Phi_0^2 \mathbb{E}(\sigma^2(X_0))D_m}{\pi_0 n \Delta} + 32c'\Delta.$$

Now, all that remains is to deal with $\Omega_n^c$. Since $\|\hat{b}_m - b_A\|_n^2 \leq \|\hat{b}_m - b\|_n^2$, it is enough to check that $\mathbb{E}(\|\hat{b}_m - b\|_n^2 \mathbb{1}_{\Omega_n^c}) \leq c/n$. Write the regression model as $Y_{k\Delta} = b(X_{k\Delta}) + \varepsilon_{k\Delta}$ with

$$\varepsilon_{k\Delta} = \frac{1}{\Delta}\int_{k\Delta}^{(k+1)\Delta}[b(X_s) - b(X_{k\Delta})]\,\mathrm{d}s + \frac{1}{\Delta}\int_{k\Delta}^{(k+1)\Delta}\sigma(X_s)\,\mathrm{d}W_s.$$

Recall that $\Pi_m$ denotes the orthogonal projection (with respect to the inner product of $\mathbb{R}^n$) onto the subspace $\{(t(X_\Delta), \ldots, t(X_{n\Delta}))', t \in S_m\}$ of $\mathbb{R}^n$. We have $(\hat{b}_m(X_\Delta), \ldots, \hat{b}_m(X_{n\Delta}))' = \Pi_m Y$, where $Y = (Y_\Delta, \ldots, Y_{n\Delta})'$. Using the same notation for the function $t$ and the vector $(t(X_\Delta), \ldots, t(X_{n\Delta}))'$, we see that

$$\|b - \hat{b}_m\|_n^2 = \|b - \Pi_m b\|_n^2 + \|\Pi_m \varepsilon\|_n^2 \leq \|b\|_n^2 + n^{-1}\sum_{i=1}^n \varepsilon_{i\Delta}^2.$$

Therefore,

$$\mathbb{E}(\|b - \hat{b}_m\|_n^2 \mathbb{1}_{\Omega_n^c}) \leq \mathbb{E}(\|b\|_n^2 \mathbb{1}_{\Omega_n^c}) + \frac{1}{n}\sum_{k=1}^n \mathbb{E}(\varepsilon_{k\Delta}^2 \mathbb{1}_{\Omega_n^c})$$

$$\leq (\mathbb{E}^{1/2}(b^4(X_0)) + \mathbb{E}^{1/2}(\varepsilon_\Delta^4))\mathbb{P}^{1/2}(\Omega_n^c).$$

By Assumption 1(ii) we have $\mathbb{E}(b^4(X_0)) \leq c(1 + \mathbb{E}(X_0^4)) = K$. With the Burholder–Davis–Gundy inequality, we find

$$\mathbb{E}(\varepsilon_\Delta^4) \leq 2^3\left\{\frac{1}{\Delta}\int_0^\Delta \mathbb{E}[(b(X_s) - b(X_\Delta))^4]\,\mathrm{d}s + \frac{36}{\Delta^3}\mathbb{E}\left(\int_0^\Delta \sigma^4(X_s)\,\mathrm{d}s\right)\right\}.$$

Under Assumptions 1, 2(i) and 3 and inequality (3), we obtain $\mathbb{E}(\varepsilon_\Delta^4) \leq C(1+\sigma_1^4/\Delta^2) := C'/\Delta^2$. The next lemma enables us to complete the proof.

**Lemma 1.** *Let $\Omega_n$ be defined by (37) and assume that $n\Delta_n/\ln^2(n) \to +\infty$ when $n \to +\infty$. Then, if $N_n \leq O(n\Delta_n/\ln^2(n))$ for collections [DP] and [W], and if $N_n \leq O((n\Delta_n)^{1/2}/\ln(n))$ for collection [T], then*

$$\mathbb{P}(\Omega_n^c) \leq \frac{c}{n^4}. \tag{39}$$

The proof of Lemma 1 is given in Section 7.

Now, we gather all terms and use (39) to obtain (15).

### 6.2. Proof of Theorem 1

The proof relies on the following Bernstein-type inequality:

**Lemma 2.** *Under the assumptions of Theorem 1, for any positive numbers $\epsilon$ and $v$, we have*

$$\mathbb{P}\left(\sum_{k=1}^n t(X_{k\Delta})Z_{k\Delta} \geq n\epsilon, \|t\|_n^2 \leq v^2\right) \leq \exp\left(-\frac{n\Delta\epsilon^2}{2\sigma_1^2 v^2}\right).$$

**Proof.** We use the fact that $\sum_{k=1}^n t(X_{k\Delta})Z_{k\Delta}$ can be written as a stochastic integral. Consider the process

$$H_u^n = H_u = \sum_{k=1}^n \mathbb{1}_{[k\Delta,(k+1)\Delta[}(u)t(X_{k\Delta})\sigma(X_u),$$

which satisfies $H_u^2 \leq \sigma_1^2 \|t\|_\infty^2$ for all $u \geq 0$. Then, writing $M_s = \int_0^s H_u \, dW_u$, we obtain that

$$M_{(n+1)\Delta} = \sum_{k=1}^n t(X_{k\Delta}) \int_{k\Delta}^{(k+1)\Delta} \sigma(X_s) \, dW_s,$$

$$\langle M \rangle_{(n+1)\Delta} = \sum_{k=1}^n t^2(X_{k\Delta}) \int_{k\Delta}^{(k+1)\Delta} \sigma^2(X_s) \, ds.$$

Moreover, $\langle M \rangle_s = \int_0^s H_u^2 \, du \leq n\sigma_1^2 \Delta \|t\|_n^2$, for all $s \geq 0$, so that $(M_s)$ and $\exp(\lambda M_s - \lambda^2 \langle M \rangle_s / 2)$ are martingales with respect to the filtration $\mathcal{F}_s = \sigma(X_u, u \leq s)$. Therefore, for all $s \geq 0$, $c > 0$, $d > 0$, $\lambda > 0$,

$$\mathbb{P}(M_s \geq c, \langle M \rangle_s \leq d) \leq \mathbb{P}\left(\exp\left(\lambda M_s - \frac{\lambda^2}{2}\langle M \rangle_s\right) \geq \exp\left(\lambda c - \frac{\lambda^2}{2}d\right)\right)$$

$$\leq \exp\left(-\left(\lambda c - \frac{\lambda^2}{2}d\right)\right).$$

Therefore,

$$\mathbb{P}(M_s \geq c, \langle M \rangle_s \leq d) \leq \inf_{\lambda > 0} \exp\left(-\left(\lambda c - \frac{\lambda^2}{2}d\right)\right) = \exp\left(-\frac{c^2}{2d}\right).$$

Finally,

$$\mathbb{P}\left(\sum_{k=1}^{n} t(X_{k\Delta})Z_{k\Delta} \geq n\epsilon, \|t\|_n^2 \leq v^2\right) = \mathbb{P}(M_{(n+1)\Delta} \geq n\Delta\epsilon, \langle M \rangle_{(n+1)\Delta} \leq nv^2\sigma_1^2\Delta)$$

$$\leq \exp\left(-\frac{(n\Delta\epsilon)^2}{2nv^2\sigma_1^2\Delta}\right) = \exp\left(-\frac{n\epsilon^2\Delta}{2v^2\sigma_1^2}\right). \quad \square$$

Now we turn to the proof of Theorem 1. As in the proof of Proposition 1, we have to split $\|\hat{b}_{\hat{m}} - b_A\|_n^2 = \|\hat{b}_{\hat{m}} - b_A\|_n^2 \mathbb{1}_{\Omega_n} + \|\hat{b}_{\hat{m}} - b_A\|_n^2 \mathbb{1}_{\Omega_n^c}$. For the treatment of $\Omega_n^c$, the end of the proof of Proposition 1 can be used.

We now focus on what happens on $\Omega_n$. From the definition of $\hat{b}_{\hat{m}}$, we have, for all $m \in \mathcal{M}_n$, $\gamma_n(\hat{b}_{\hat{m}}) + \text{pen}(\hat{m}) \leq \gamma_n(b_m) + \text{pen}(m)$. We proceed as in the proof of Proposition 1 with the additional penalty terms (see (38)) and obtain

$$\mathbb{E}(\|\hat{b}_{\hat{m}} - b_A\|_n^2 \mathbb{1}_{\Omega_n}) \leq 7\pi_1 \|b_m - b_A\|^2 + 4\,\text{pen}(m) + 32\mathbb{E}\left(\sup_{t \in S_m + S_{\hat{m}}, \|t\|_\pi = 1} [\nu_n(t)]^2 \mathbb{1}_{\Omega_n}\right)$$

$$- 4\mathbb{E}(\text{pen}(\hat{m})) + 32c'\Delta.$$

The main problem here is to control the supremum of $\nu_n(t)$ on a random ball (which depends on the random $\hat{m}$). This is done by using the martingale property of $\nu_n(t)$.

Let us introduce the notation

$$G_m(m') = \sup_{t \in S_m + S_{m'}, \|t\|_\pi = 1} |\nu_n(t)|.$$

Now, we plug in a function $p(m, m')$, which will in turn fix the penalty:

$$G_m^2(\hat{m})\mathbb{1}_{\Omega_n} \leq [(G_m^2(\hat{m}) - p(m, \hat{m}))\mathbb{1}_{\Omega_n}]_+ + p(m, \hat{m})$$

$$\leq \sum_{m' \in \mathcal{M}_n} [(G_m^2(m') - p(m, m'))\mathbb{1}_{\Omega_n}]_+ + p(m, \hat{m}).$$

And pen is chosen such that $8p(m, m') \leq \text{pen}(m) + \text{pen}(m')$. More precisely, the next proposition determines the choice of $p(m, m')$.

**Proposition 3.** *Under the assumptions of Theorem 1, there exists a numerical constant $\kappa_1$ such that, for $p(m,m') = \kappa_1 \sigma_1^2 (D_m + D_{m'})/(n\Delta)$, we have*

$$\mathbb{E}[(G_m^2(m') - p(m,m'))\mathbb{1}_{\Omega_n}]_+ \leq c\sigma_1^2 \frac{e^{-D_{m'}}}{n\Delta}.$$

**Proof of Proposition 3.** The result of Proposition 3 follows from the inequality of Lemma 2 by the $\mathbb{L}^2$-chaining technique used in Baraud *et al.* [5] (see their Section 7, pages 44–47, Lemma 7.1, with $s^2 = \sigma_1^2/\Delta$). □

It is easy to see that the result of Theorem 1 follows from Proposition 3 with $\text{pen}(m) \geq \kappa \sigma_1^2 D_m/(n\Delta)$ and $\kappa = 8\kappa_1$.

### 6.3. Proof of Proposition 2

First, we prove that

$$\mathbb{E}\left(\frac{1}{n}\sum_{k=1}^{n} R_{k\Delta}^2\right) \leq K\Delta^2. \tag{40}$$

With the obvious convention, let $R_{k\Delta} = R_{k\Delta}^{(1)} + R_{k\Delta}^{(2)} + R_{k\Delta}^{(3)}$ so that (40) holds if $\mathbb{E}[(R_{k\Delta}^{(i)})^2] \leq K_i \Delta^2$ for $i = 1, 2, 3$. Using Assumption 1,

$$\mathbb{E}[(R_{k\Delta}^{(1)})^2] \leq \mathbb{E}\left(\int_{k\Delta}^{(k+1)\Delta} b^2(X_s)\,ds\right)^2 \leq \Delta \mathbb{E}\left(\int_{k\Delta}^{(k+1)\Delta} b^4(X_s)\,ds\right)$$
$$\leq \Delta^2 \mathbb{E}(b^4(X_0)) \leq c\Delta^2.$$

We also have

$$\mathbb{E}[(R_{k\Delta}^{(2)})^2] \leq \frac{1}{\Delta^2}\left(\mathbb{E}\left(\int_{k\Delta}^{(k+1)\Delta}(b(X_s) - b(X_{k\Delta}))\,ds\right)^4 \mathbb{E}\left(\int_{k\Delta}^{(k+1)\Delta} \sigma(X_s)\,dW_s\right)^4\right)^{1/2}.$$

Using (3), we obtain

$$\mathbb{E}[(R_{k\Delta}^{(2)})^2] \leq c'\Delta^2.$$

Lastly, using Assumptions 1 and 2 and equation (22),

$$\mathbb{E}[(R_{k\Delta}^{(3)})^2] \leq \frac{1}{\Delta}\mathbb{E}\left(\int_{k\Delta}^{(k+1)\Delta}((k+1)\Delta - s)^2 \psi^2(X_s)\,ds\right) \leq \mathbb{E}(\psi^2(X_0))\frac{\Delta^2}{3} \leq c''\Delta^2.$$

Therefore (40) is proved.

We now return to (23) and recall that $\Omega_n$ is defined by (37). The treatment is similar to that for the drift estimator. On $\Omega_n$, $\|\hat{\sigma}_m^2 - \sigma_m^2\|_\pi^2 \leq 2\|\hat{\sigma}_m^2 - \sigma_m^2\|_n^2$,

$$\|\hat{\sigma}_m^2 - \sigma_A^2\|_n^2 \leq \|\sigma_m^2 - \sigma_A^2\|_n^2 + \frac{1}{8}\|\hat{\sigma}_m^2 - \sigma_m^2\|_\pi^2 + 8 \sup_{t \in S_m, \|t\|_\pi = 1} \breve{\nu}_n^2(t)$$

$$+ \frac{1}{8}\|\hat{\sigma}_m^2 - \hat{\sigma}_m^2\|_n^2 + \frac{8}{n}\sum_{k=1}^n R_{k\Delta}^2$$

$$\leq \|\sigma_m^2 - \sigma_A^2\|_n^2 + \frac{3}{8}\|\hat{\sigma}_m^2 - \sigma_m^2\|_n^2 + 8 \sup_{t \in S_m, \|t\|_\pi = 1} \breve{\nu}_n^2(t) + \frac{8}{n}\sum_{k=1}^n R_{k\Delta}^2.$$

Setting $B_m(0,1) = \{t \in S_m, \|t\| \leq 1\}$ and $B_m^\pi(0,1) = \{t \in S_m, \|t\|_\pi \leq 1\}$, the following holds on $\Omega_n$:

$$\frac{1}{4}\|\hat{\sigma}_m^2 - \sigma_A^2\|_n^2 \leq \frac{7}{4}\|\sigma_m^2 - \sigma_A^2\|_n^2 + 8 \sup_{t \in B_m^\pi(0,1)} \breve{\nu}_n^2(t) + \frac{8}{n}\sum_{k=1}^n R_{k\Delta}^2.$$

Moreover,

$$\mathbb{E}\left(\sup_{t \in B_m^\pi(0,1)} \breve{\nu}_n^2(t)\right) \leq \frac{1}{\pi_0}\mathbb{E}\left(\sup_{t \in B_m(0,1)} \breve{\nu}_n^2(t)\right) \leq \frac{1}{\pi_0}\sum_{\lambda \in \Lambda_m} \mathbb{E}(\breve{\nu}_n^2(\varphi_\lambda))$$

$$\leq \frac{1}{\pi_0 n^2} \sum_{\lambda \in \Lambda_m} \mathbb{E}\left(\sum_{k=1}^n \varphi_\lambda^2(X_{k\Delta}) V_{k\Delta}^2\right)$$

$$\leq \frac{\Phi_0^2 D_m}{\pi_0 n}\{12\mathbb{E}(\sigma^4(X_0)) + 4\Delta C_{b,\sigma}\},$$

where $C_{b,\sigma} = \mathbb{E}((\sigma'\sigma^2)^2(X_0)) + \sigma_1^2\mathbb{E}(b^2(X_0))$. Now using the condition on $N_n$, we have $\Delta D_m/n \leq \Delta N_n/n \leq \Delta^2/\ln^2(n)$. This yields the first three terms of the right-hand side of (24).

The treatment of $\Omega_n^c$ is the same as for $b$ with the regression model $U_{k\Delta} = \sigma^2(X_{k\Delta}) + \eta_{k\Delta}$, where $\eta_{k\Delta} = V_{k\Delta} + R_{k\Delta}$. By standard inequalities, $\mathbb{E}(\eta_\Delta^4) \leq K\{\Delta^4\mathbb{E}(b^8(X_0)) + \mathbb{E}(\sigma^8(X_0))\}$. Hence, $\mathbb{E}(\eta_\Delta^4)$ is bounded. Moreover, using Lemma 1, $\mathbb{P}(\Omega_n^c) \leq c/n^2$.

### 6.4. Proof of Theorem 2

This proof follows the same lines as the proof of Theorem 1. We start with a Bernstein-type inequality.

**Lemma 3.** *Under the assumptions of Theorem 2,*

$$\mathbb{P}\left(\sum_{k=1}^n t(X_{k\Delta})V_{k\Delta}^{(1)} \geq n\epsilon, \|t\|_n^2 \leq v^2\right) \leq \exp\left(-Cn\frac{\epsilon^2/2}{2\sigma_1^4 v^2 + \epsilon\|t\|_\infty \sigma_1^2 v}\right)$$

*and*

$$\mathbb{P}\left(\frac{1}{n}\sum_{k=1}^{n} t(X_{k\Delta})V_{k\Delta}^{(1)} \geq v\sigma_1^2(2x)^{1/2} + \sigma_1^2\|t\|_\infty x, \|t\|_n^2 \leq v^2\right) \leq \exp(-Cnx). \qquad (41)$$

The non-trivial link between the above two inequalities is enhanced by Birgé and Massart [14], so we just prove the first.

**Proof of Lemma 3.** First we note that

$$\mathbb{E}(e^{ut(X_{n\Delta})V_{n\Delta}^{(1)}}|\mathcal{F}_{n\Delta}) = 1 + \sum_{p=2}^{+\infty}\frac{u^p}{p!}\mathbb{E}\{(t(X_{n\Delta})V_{n\Delta}^{(1)})^p|\mathcal{F}_{n\Delta}\}$$

$$\leq 1 + \sum_{p=2}^{+\infty}\frac{u^p}{p!}|t(X_{n\Delta})|^p\mathbb{E}(|V_{n\Delta}^{(1)}|^p|\mathcal{F}_{n\Delta}).$$

Next we apply successively the Hölder inequality and the Burkholder–Davis–Gundy inequality with best constant (Proposition 4.2 of Barlow and Yor [6]). For a continuous martingale $(M_t)$, with $M_0 = 0$, for $k \geq 2$, $M_t^* = \sup_{s \leq t}|M_s|$ satisfies $\|M^*\|_k \leq ck^{1/2}\|\langle M\rangle^{1/2}\|_k$, with $c$ a universal constant. And we obtain

$$\mathbb{E}(|V_{n\Delta}^{(1)}|^p|\mathcal{F}_{n\Delta}) \leq \frac{2^{p-1}}{\Delta^p}\left\{\mathbb{E}\left(\left|\int_{n\Delta}^{(n+1)\Delta}\sigma(X_s)\,\mathrm{d}W_s\right|^{2p}\bigg|\mathcal{F}_{n\Delta}\right)\right.$$

$$\left.+ \mathbb{E}\left(\left|\int_{n\Delta}^{(n+1)\Delta}\sigma^2(X_s)\,\mathrm{d}s\right|^p\bigg|\mathcal{F}_{n\Delta}\right)\right\}$$

$$\leq \frac{2^{p-1}}{\Delta^p}(c^{2p}(2p)^p\Delta^p\sigma_1^{2p} + \Delta^p\sigma_1^{2p}) \leq (2\sigma_1 c)^{2p}p^p.$$

Therefore,

$$\mathbb{E}(e^{ut(X_{n\Delta})V_{n\Delta}^{(1)}}|\mathcal{F}_{n\Delta}) \leq 1 + \sum_{k=2}^{\infty}\frac{p^p}{p!}(4u\sigma_1^2c^2)^p|t(X_{n\Delta})|^p.$$

Using $p^p/p! \leq e^{p-1}$, we find

$$\mathbb{E}(e^{ut(X_{n\Delta})V_{n\Delta}^{(1)}}|\mathcal{F}_{n\Delta}) \leq 1 + e^{-1}\sum_{k=2}^{\infty}(4u\sigma_1^2c^2e)^p|t(X_{n\Delta})|^p$$

$$\leq 1 + e^{-1}\frac{(4u\sigma_1^2c^2e)^2t^2(X_{n\Delta})}{1 - (4u\sigma_1^2c^2e\|t\|_\infty)}.$$

Now, let us set

$$a = e(4\sigma_1^2c^2)^2 \quad \text{and} \quad b = 4\sigma_1^2c^2e\|t\|_\infty.$$

Since for $x \geq 0$, $1 + x \leq e^x$, we obtain, for all $u$ such that $bu < 1$,

$$\mathbb{E}(e^{ut(X_{n\Delta})V^{(1)}_{n\Delta}}|\mathcal{F}_{n\Delta}) \leq 1 + \frac{au^2t^2(X_{n\Delta})}{1-bu} \leq \exp\left(\frac{au^2t^2(X_{n\Delta})}{1-bu}\right).$$

This can also be written as

$$\mathbb{E}\left(\exp\left(ut(X_{n\Delta})V^{(1)}_{n\Delta} - \frac{au^2t^2(X_{n\Delta})}{1-bu}\right)\Big|\mathcal{F}_{n\Delta}\right) \leq 1.$$

Therefore, iterating conditional expectations yields

$$\mathbb{E}\left[\exp\left\{\sum_{k=1}^{n}\left(ut(X_{k\Delta})V^{(1)}_{k\Delta} - \frac{au^2t^2(X_{k\Delta})}{1-bu}\right)\right\}\right] \leq 1.$$

Then we deduce that

$$\mathbb{P}\left(\sum_{k=1}^{n} t(X_{k\Delta})V^{(1)}_{k\Delta} \geq n\epsilon, \|t\|_n^2 \leq v^2\right)$$

$$\leq e^{-nu\epsilon}\mathbb{E}\left\{\mathbb{1}_{\|t\|_n^2 \leq v^2} \exp\left(u\sum_{k=1}^{n} t(X_{k\Delta})V^{(1)}_{k\Delta}\right)\right\}$$

$$\leq e^{-nu\epsilon}\mathbb{E}\left[\mathbb{1}_{\|t\|_n^2 \leq v^2} \exp\left\{\sum_{k=1}^{n}\left(ut(X_{k\Delta})V^{(1)}_{k\Delta} - \frac{au^2t^2(X_{k\Delta})}{1-bu}\right)\right\} e^{(au^2)/(1-bu)\sum_{k=1}^{n} t^2(X_{k\Delta})}\right]$$

$$\leq e^{-nu\epsilon} e^{(nau^2v^2)/(1-bu)}\mathbb{E}\left[\exp\left\{\sum_{k=1}^{n}\left(ut(X_{k\Delta})V^{(1)}_{k\Delta} - \frac{au^2t^2(X_{k\Delta})}{1-bu}\right)\right\}\right]$$

$$\leq e^{-nu\epsilon} e^{(nau^2v^2)/(1-bu)}.$$

The inequality holds for any $u$ such that $bu < 1$. In particular, $u = \epsilon/(2av^2 + \epsilon b)$ gives $-u\epsilon + av^2u^2/(1-bu) = -(1/2)(\epsilon^2/(2av^2 + \epsilon b))$ and therefore

$$\mathbb{P}\left(\sum_{k=1}^{n} t(X_{k\Delta})V^{(1)}_{k\Delta} \geq n\epsilon, \|t\|_n^2 \leq v^2\right) \leq \exp\left(-n\frac{\epsilon^2/2}{2av^2 + \epsilon b}\right). \qquad \square$$

As for $\hat{b}_{\hat{m}}$, we introduce the additional penalty terms and obtain that the risk satisfies

$$\mathbb{E}(\|\hat{\sigma}^2_{\hat{m}} - \sigma^2_A\|_n^2 \mathbb{1}_{\Omega_n}) \leq 7\pi_1\|\sigma^2_m - \sigma^2_A\|^2 + 4\widetilde{\mathrm{pen}}(m) + 32\mathbb{E}\left(\sup_{t \in B^{\pi}_{m,\hat{m}}(0,1)} (\check{\nu}_n(t))^2 \mathbb{1}_{\Omega_n}\right)$$

$$- 4\mathbb{E}(\widetilde{\mathrm{pen}}(\hat{m})) + K'\Delta^2, \tag{42}$$

where $B_{m,m'}^\pi(0,1) = \{t \in S_m + S_{m'}, \|t\|_\pi = 1\}$. Let us denote by

$$\breve{G}_m(m') = \sup_{t \in B_{m,m'}^\pi(0,1)} |\breve{\nu}_n^{(1)}(t)|$$

the main quantity to be studied, where

$$\breve{\nu}_n^{(1)}(t) = \frac{1}{n} \sum_{k=1}^n t(X_{k\Delta}) V_{k\Delta}^{(1)};$$

define also

$$\breve{\nu}_n^{(2)}(t) = \frac{1}{n} \sum_{k=1}^n t(X_{k\Delta})(V_{k\Delta}^{(2)} + V_{k\Delta}^{(3)}).$$

As for the drift, we write

$$\mathbb{E}(\breve{G}_m^2(\hat{m})) \leq \mathbb{E}[(\breve{G}_m^2(\hat{m}) - \tilde{p}(m,\hat{m}))\mathbb{1}_{\Omega_n}]_+ + \mathbb{E}(\tilde{p}(m,\hat{m}))$$
$$\leq \sum_{m' \in \mathcal{M}_n} \mathbb{E}[(\breve{G}_m^2(m') - \tilde{p}(m,m'))\mathbb{1}_{\Omega_n}]_+ + \mathbb{E}(\tilde{p}(m,\hat{m})).$$

Now we have the following statement.

**Proposition 4.** *Under the assumptions of Theorem 2, for*

$$\tilde{p}(m,m') = \kappa^* \sigma_1^4 \left\{ \frac{D_m + D_{m'}}{n} + \frac{\Phi_0^2}{\pi_0} \left( \frac{D_m + D_{m'}}{n} \right)^2 \right\},$$

*where $\kappa^*$ is a numerical constant, we have*

$$\mathbb{E}[(\breve{G}_m^2(m') - \tilde{p}(m,m'))\mathbb{1}_{\Omega_n}]_+ \leq c\sigma_1^4 \frac{\mathrm{e}^{-D_{m'}}}{n}.$$

The result of Proposition 4 is obtained from inequality (41) of Lemma 3 by a $L^2(\pi) - L^\infty$ chaining technique. For a description of this method, in a more general setting, we refer to Propositions 2–4 in Comte ([15], page 282–287), to Theorem 5 in Birgé and Massart [14] and to Proposition 7 and Theorems 8 and 9 in Barron *et al.* [8]. Note that there is a difference between Propositions 3 and 4 which comes from the additional term $\|t\|_\infty$ appearing in Lemma 3. For this reason, we need to use the fact that $\|\sum_{\lambda \in \Lambda_m} \beta_\lambda \psi_\lambda\|_\infty / \sup_{\lambda \in \Lambda_m} |\beta_\lambda| \leq \|\sum |\psi_\lambda|\|_\infty \leq (r_{\max}+1) D_m^{1/2}/\pi_0^{1/2}$ for $(\psi_\lambda)_{\lambda \in \Lambda_m}$ an $\mathbb{L}^2(\pi)$-orthonormal basis constructed by orthonormalisation of the $(\varphi_\lambda)$. This explains the additional term appearing in $\tilde{p}(m,m')$.

Choosing $\widetilde{\mathrm{pen}}(m) \geq \tilde{\kappa} \sigma_1^4 D_m/n$ with $\tilde{\kappa} = 16\kappa^*$, we deduce from (42), Proposition 4 and $D_m \leq N_n \leq n\Delta/\ln^2(n)$ that

$$\mathbb{E}(\|\hat{\sigma}_{\hat{m}}^2 - \sigma_A^2\|_n^2) \leq 7\pi_1 \|\sigma_m^2 - \sigma_A^2\|^2 + 8\widetilde{\mathrm{pen}}(m) + c\sigma_1^4 \sum_{m' \in \mathcal{M}_n} \frac{\mathrm{e}^{-D_{m'}}}{n} + \frac{\tilde{\kappa}\sigma_1^4 \Phi_0^2}{\pi_0} \frac{\Delta^2}{\ln^4(n)}$$

$$+ 64\mathbb{E}\left(\sup_{t\in B^\pi_{m,\hat{m}}(0,1)} (\breve{\nu}_n^{(2)}(t))^2\right) + K'\Delta^2 + \mathbb{E}(\|\hat{\sigma}_{\hat{m}}^2 - \sigma_A^2\|_n^2 \mathbb{1}_{\Omega_n^c}).$$

The bound for $\mathbb{E}(\|\hat{\sigma}_{\hat{m}}^2 - \sigma^2\|_n^2 \mathbb{1}_{\Omega_n^c})$ is the same as that given in the end of the proof of Proposition 2. It is less than $c/n$ provided that $N_n \leq n\Delta/\ln^2(n)$ for [DP] and [W] and $N_n^2 \leq n\Delta/\ln^2(n)$ for [T].

Since the spaces are all contained in a space denoted by $\mathcal{S}_n$ with dimension $N_n$ bounded as right above, we have

$$\mathbb{E}\left(\sup_{t\in B^\pi_{m,\hat{m}}(0,1)} (\breve{\nu}_n^{(2)}(t))^2\right) \leq \frac{1}{\pi_0}\mathbb{E}\left(\sup_{t\in\mathcal{S}_n, \|t\|=1} (\breve{\nu}_n^{(2)}(t))^2\right) \leq KC_{b,\sigma}\Phi_0^2 \frac{\Delta N_n}{\pi_0 n} \leq K'\Delta^2.$$

The result of Theorem 2 follows.

## 7. Proof of Lemma 1

Using Baraud *et al.* [4], we prove that, for all $n$ and $\Delta > 0$,

$$\mathbb{P}(\Omega_n^c) \leq 2n\beta_X(q_n\Delta) + 2n^2 \exp\left(-C_0 \frac{n}{q_n L_n(\phi)}\right), \tag{43}$$

where $C_0$ is a constant depending on $\pi_0, \pi_1$, $q_n$ is an integer such that $q_n < n$, and $L_n(\phi)$ is a quantity depending on the basis of the largest nesting space $\mathcal{S}_n$ of the collection and is defined below. We recall that $N_n = \dim(\mathcal{S}_n)$.

We first prove (43). We use Berbee's coupling method as in Proposition 5.1 of Viennet [32] and its proof. We assume that $n = 2p_n q_n$. Then there exist random variables $X^*_{i\Delta}$, $i = 1, \ldots, n$, satisfying the following properties:

- For $\ell = 1, \ldots, p_n$, the random vectors $\vec{U}_{\ell,1} = (X_{[2(\ell-1)q_n+1]\Delta}, \ldots, X_{(2\ell-1)q_n\Delta})'$ and $\vec{U}^*_{\ell,1} = (X^*_{[2(\ell-1)q_n+1]\Delta}, \ldots, X^*_{(2\ell-1)q_n\Delta})'$ have the same distribution, and so have the vectors $\vec{U}_{\ell,2} = (X_{[(2\ell-1)q_n+1]\Delta}, \ldots, X_{2\ell q_n\Delta})'$ and $\vec{U}^*_{\ell,2} = (X^*_{[(2\ell-1)q_n+1]\Delta}, \ldots, X^*_{2\ell q_n\Delta})'$.
- For $\ell = 1, \ldots, p_n$, $\mathbb{P}(\vec{U}_{\ell,1} \neq \vec{U}^*_{\ell,1}) \leq \beta_X(q_n\Delta)$ and $\mathbb{P}(\vec{U}_{\ell,2} \neq \vec{U}^*_{\ell,2}) \leq \beta_X(q_n\Delta)$.
- For each $\delta \in \{1,2\}$, the random vectors $\vec{U}^*_{1,\delta}, \ldots, \vec{U}^*_{p_n,\delta}$ are independent.

Let us define $\Omega^* = \{X_{i\Delta} = X^*_{i\Delta}, i = 1, \ldots, n\}$. We have $\mathbb{P}(\Omega_n^c) \leq \mathbb{P}(\Omega_n^c \cap \Omega^*) + \mathbb{P}(\Omega^{*c})$ and clearly

$$\mathbb{P}(\Omega^{*c}) \leq 2p_n\beta_X(q_n\Delta) \leq n\beta_X(q_n\Delta). \tag{44}$$

Thus, (43) holds if we prove

$$\mathbb{P}(\Omega_n^c \cap \Omega^*) \leq 2N_n^2 \exp\left(-A_0 \frac{\pi_0^2}{\pi_1} \frac{n}{q_n L_n(\phi)}\right), \tag{45}$$

where $L_n(\phi)$ is defined as follows. Let $(\varphi_\lambda)_{\lambda \in \Lambda_n}$ be an $\mathbb{L}^2(A, \mathrm{d}x)$-orthonormal basis of $\mathcal{S}_n$ and, as in Baraud *et al.* [4], define the matrices

$$V = \left[\left(\int_A \varphi_\lambda^2(x)\varphi_{\lambda'}^2(x)\,\mathrm{d}x\right)^{1/2}\right]_{\lambda,\lambda' \in \Lambda_n \times \Lambda_n}, \qquad B = (\|\varphi_\lambda \varphi_{\lambda'}\|_\infty)_{\lambda,\lambda' \in \Lambda_n \times \Lambda_n}.$$

Then we set $L_n(\phi) = \max\{\rho^2(V), \rho(B)\}$, where, for any symmetric matrix $M = (M_{\lambda,\lambda'})$, $\rho(M) = \sup_{\{a_\lambda\}, \sum_\lambda a_\lambda^2 \leq 1} \sum_{\lambda,\lambda'} |a_\lambda||a_{\lambda'}||M_{\lambda,\lambda'}|$.

We now prove (45). Let $\mathbb{P}^*(\cdot) := \mathbb{P}(\cdot \cap \Omega^*)$. We use Baraud ([3], Claim 2 in Proposition 4.2). Consider $v_n(t) = (1/n)\sum_{i=1}^n [t(X_{i\Delta}) - \mathbb{E}(t(X_{i\Delta}))]$, $B_\pi(0,1) = \{t \in \mathcal{S}_n, \|t\|_\pi \leq 1\}$ and $B(0,1) = \{t \in \mathcal{S}_n, \|t\| \leq 1\}$. As, on $A$, $\pi_0 \leq \pi(x) \leq \pi_1$,

$$\sup_{t \in B_\pi(0,1)} |v_n(t^2)| = \sup_{t \in \mathcal{S}_n / \{0\}} \left|\frac{\|t\|_n^2}{\|t\|_\pi^2} - 1\right| \leq \pi_0^{-1} \sup_{t \in B(0,1)} |v_n(t^2)|.$$

Thus

$$\mathbb{P}^*\left(\sup_{t \in B_\pi(0,1)} |v_n(t^2)| \geq \rho_0\right) \leq \mathbb{P}^*\left(\sup_{t \in B(0,1)} |v_n(t^2)| \geq \pi_0 \rho_0\right)$$

$$\leq \mathbb{P}^*\left(\sup_{\sum_{\lambda \in \Lambda_n} a_\lambda^2 \leq 1} \sum_{\lambda,\lambda' \in \Lambda_n} |a_\lambda a_{\lambda'}||v_n(\varphi_\lambda \varphi_{\lambda'})| \geq \pi_0 \rho_0\right).$$

On the set $\{\forall (\lambda, \lambda') \in \Lambda_n^2, |v_n(\varphi_\lambda \varphi_{\lambda'})| \leq 2V_{\lambda\lambda'}(2\pi_1 x)^{1/2} + 3B_{\lambda\lambda'} x\}$, we have

$$\sup_{\sum_{\lambda \in \Lambda_n} a_\lambda^2 \leq 1} \sum_{\lambda,\lambda' \in \Lambda_n} |a_\lambda a_{\lambda'}||v_n(\varphi_\lambda \varphi_{\lambda'})| \leq 2\rho(V)(2\pi_1 x)^{1/2} + 3\rho(B)x.$$

By choosing $x = (\rho_0 \pi_0)^2/(16\pi_1 L_n(\phi))$ and $\rho_0 = 1/2$, and recall that $\pi_0 \leq \pi_1$, we obtain that

$$\sup_{\sum_{\lambda \in \Lambda_n} a_\lambda^2 \leq 1} \sum_{\lambda\lambda'} |a_\lambda a_{\lambda'}||v_n(\varphi_\lambda \varphi_{\lambda'})| \leq \rho_0 \pi_0 = \frac{\pi_0}{2}.$$

This leads to

$$\mathbb{P}^*(\Omega_n^c) = \mathbb{P}^*\left(\sup_{t \in B_\pi(0,1)} |v_n(t^2)| \geq \frac{1}{2}\right)$$

$$\leq \mathbb{P}^*(\{\forall (\lambda, \lambda') \in \Lambda_n^2, |v_n(\varphi_\lambda \varphi_{\lambda'})| \geq 2V_{\lambda\lambda'}(2\pi_1 x)^{1/2} + 3B_{\lambda\lambda'} x\}).$$

The proof of (45) is then achieved by using the following claim, which is exactly Claim 6 in the proof of Proposition 7 of Baraud *et al.* [4].

**Claim 1.** *Let $(\varphi_\lambda)_{\lambda \in \Lambda_n}$ be an $\mathbb{L}^2(A, \mathrm{d}x)$ orthonormal basis of $\mathcal{S}_n$. Then, for all $x \geq 0$ and all integers $q$, $1 \leq q \leq n$,*

$$\mathbb{P}^*(\exists (\lambda, \lambda') \in \Lambda_n^2 / |v_n(\varphi_\lambda \varphi_{\lambda'})| > 2V_{\lambda, \lambda'}(2\pi_1 x)^{1/2} + 2B_{\lambda, \lambda'} x) \leq 2N_n^2 \exp\left(-\frac{nx}{q}\right).$$

Claim 1 implies that

$$\mathbb{P}(\Omega_n^c \cap \Omega^*) \leq 2N_n^2 \exp\left(-\frac{\pi_0^2}{64\pi_1} \frac{n}{q_n L_n(\phi)}\right),$$

and thus (45) holds true.

Again we refer to Baraud *et al.* [4] (see Lemma 2 in Section 10). It is proved there that, for [T], $L_n(\phi) \leq C_\phi N_n^2$. For [W] and [DP] (see Sections 2.2 and 2.3 above), $L_n(\phi) \leq C'_\phi N_n$. We now use (43) to complete the proof of Lemma 1. By assumption, the diffusion process $X$ is geometrically $\beta$-mixing. So, for some constant $\theta$, $\beta_X(q_n\Delta) \leq \mathrm{e}^{-\theta q_n \Delta}$. Provided that $\Delta = \Delta_n$ satisfies $\ln(n)/(n\Delta) \to 0$, it is possible to take $q_n = [5\ln(n)/(\theta\Delta)] + 1$. This yields

$$\mathbb{P}(\Omega_n^c) \leq \frac{2}{n^4} + 2n^2 \exp\left(-C'_0 \frac{n\Delta}{\ln(n) N_n}\right).$$

The above constraint on $\Delta$ must be strengthened. Indeed, to ensure (39), we need

$$\frac{n\Delta}{N_n} \geq \frac{6\ln^2(n)}{C'_0}, \quad \text{i.e.} \quad N_n \leq \tilde{C}_0 \frac{n\Delta}{\ln^2(n)}$$

for [W] and [DP]. This requires $n\Delta/\ln^2(n) \to +\infty$. The result for [T] follows analogously.

## Acknowledgements

The authors wish to thank the Associate Editor and anonymous referees for comments that helped to significantly improve the paper.